\documentclass[a4paper,12pt]{amsart}

\usepackage[english,francais]{babel}

\usepackage[utf8]{inputenc} 

\usepackage[T1]{fontenc}

\usepackage{amsmath}
\usepackage{amsfonts}
\usepackage{amssymb}

\usepackage{xy}
\usepackage{bm}
\xyoption{all}
\xyoption{poly}
\xyoption{arc}

\usepackage{epsfig}
\usepackage{a4wide}

\usepackage{url} 
\usepackage{enumerate} 

\usepackage{hyperref}

\newtheorem{thm}{Théorème}[section]
\newtheorem{cor}[thm]{Corollaire}
\newtheorem{lemme}[thm]{Lemme}
\newtheorem{prop}[thm]{Proposition}
\theoremstyle{definition}

\newtheorem{defn}[thm]{Définition}
\newtheorem{rqe}[thm]{Remarque}

\newcommand{\Vreg}{{V^{\mathrm{reg}}}}

\newcommand{\Enreg}{{E_n^{\mathrm{reg}}}}

\newcommand{\eh}{{e^{2i\pi/h}}}
\newcommand{\tq}{{~|~}}
\newcommand{\ie}{{\emph{i.e.}~}}
\newcommand{\surj}{{\twoheadrightarrow}}
\newcommand{\inj}{{\hookrightarrow}}
\newcommand{\ssi}{{\Leftrightarrow}}

\newcommand{\FS}{{\mathfrak{S}}} 
 
\newcommand{\CO}{{\mathcal{O}}} 
\newcommand{\CL}{{\mathcal{L}}} 
 \newcommand{\CH}{{\mathcal{H}}}
 \newcommand{\CR}{{\mathcal{R}}}
 
\newcommand{\CU}{{\mathcal{U}}} \newcommand{\CA}{{\mathcal{A}}}

\renewcommand{\CR}{{\mathcal{R}}}
\renewcommand{\CL}{{\mathcal{L}}} 

\newcommand{\CK}{{\mathcal{K}}}
 
\renewcommand{\CU}{{\mathcal{U}}} 
\newcommand{\CB}{{\mathcal{B}}}

\newcommand{\BN}{{\mathbb{N}}}
\newcommand{\BC}{{\mathbb{C}}} 
\newcommand{\BZ}{{\mathbb{Z}}}

\newcommand{\qg}{{\backslash}} % "quotient à gauche"
\newcommand{\eps}{\varepsilon} 
\newcommand{\D}{\Delta}
\newcommand{\fconj}{\stackrel{c}{\sim}}
\newcommand{\Lb}{\bar{\CL}}

\newcommand{\GL}{\operatorname{GL}\nolimits}

\newcommand{\OO}{\operatorname{O}\nolimits}
\newcommand{\U}{\operatorname{U}\nolimits}

\newcommand{\Red}{\operatorname{Red}\nolimits}

\newcommand{\Cat}{\operatorname{Cat}\nolimits}
\newcommand{\reg}{\operatorname{reg}\nolimits}

\newcommand{\im}{\operatorname{im}\nolimits}
\newcommand{\re}{\operatorname{re}\nolimits}
\newcommand{\Spec}{\operatorname{Spec}\nolimits}
\newcommand{\Ker}{\operatorname{Ker}\nolimits}

\renewcommand{\Im}{\operatorname{Im}\nolimits}

\newcommand{\codim}{\operatorname{codim}\nolimits}

\newcommand{\Disc}{\operatorname{Disc}\nolimits}

\renewcommand{\Red}{\operatorname{Red}\nolimits}

\newcommand{\Fact}{\operatorname{\textsc{fact}}\nolimits}
\newcommand{\fact}{\operatorname{fact}\nolimits}
\newcommand{\LL}{\operatorname{LL}\nolimits}
\newcommand{\lbl}{\operatorname{lbl}\nolimits}

\newcommand{\rg}{\operatorname{rg}\nolimits}
\newcommand{\NCP}{\operatorname{\textsc{ncp}}\nolimits}
\newcommand{\cp}{\operatorname{cp}\nolimits}
\newcommand{\comp}{\operatorname{comp}\nolimits}

\newcommand{\<}{\operatorname{\preccurlyeq}}

\newcommand{\lex}{\mathrel{\leq_{\mathrm{lex}}}}
\newcommand{\pr}{\operatorname{pr}\nolimits}

\begin{document}

\title{Orbites d'Hurwitz des factorisations primitives d'un élément de Coxeter}

\author{Vivien Ripoll}
%Institution sur le nom ?
 \address{DMA, École normale supérieure, 
   45 rue d'Ulm, 75230 Paris cedex 05, France}
\email{vivien.ripoll@ens.fr}

\begin{abstract}
  On considère l'action d'Hurwitz du groupe de tresses
  usuel sur les factorisations d'un élément de Coxeter $c$ d'un groupe de
  réflexions complexe bien engendré $W$. Il est bien connu que l'action
  d'Hurwitz est transitive sur l'ensemble des décompositions réduites de
  $c$ en réflexions. On démontre ici une propriété similaire pour les
  factorisations primitives de $c$, \ie celles dont tous les facteurs sauf
  un sont des réflexions. Cette étude est motivée par la recherche d'une
  explication géométrique de la formule de Chapoton sur le nombre de
  chaînes de longueur donnée dans le treillis des partitions non croisées
  $\NCP_W$. La démonstration présentée repose sur les propriétés du
  revêtement de Lyashko-Looijenga et sur la géométrie du discriminant de
  $W$.

  \textbf{Mots-clés :} action d'Hurwitz, groupe de réflexions complexe, élément
  de Coxeter, treillis des partitions non croisées, revêtement de
  Lyashko-Looijenga.
  
  \medskip
  \selectlanguage{english}
  \noindent \textsc{Abstract}. We study the Hurwitz action of the classical
  braid group on factorisations of a Coxeter element $c$ in a
  well-generated complex reflection group $W$. It is well known that the
  Hurwitz action is transitive on the set of reduced decompositions of $c$
  in reflections. Our main result is a similar property for the primitive
  factorisations of $c$, \ie factorisations with only one factor which is
  not a reflection. The motivation is the search for a geometric proof of
  Chapoton's formula for the number of chains of given length in the
  non-crossing partitions lattice $\NCP_W$. Our proof uses the properties
  of the Lyashko-Looijenga covering and the geometry of the discriminant of
  $W$.

 \selectlanguage{francais}
\end{abstract}

\maketitle

\section*{Introduction}

Soit $V$ un $\BC$-espace vectoriel de dimension finie, et $W \subseteq
\GL(V)$ un groupe de réflexions complexe fini bien engendré (définition
en partie \ref{subpartcox}). On note $\CR$ l'ensemble de
toutes les réflexions de $W$.  On définit la longueur d'un élément $w$ de $W$ :
\[ \ell(w) := \min \{k \in \BN \tq \exists r_1,\dots, r_k \in \CR,\
w=r_1\dots r_k \} .\]
On note $\Red (w)$ l'ensemble des décompositions réduites de $w$, \ie des
mots de longueur $\ell(w)$ représentant $w$. Ce sont donc des factorisations
minimales de $w$ en réflexions. Plus généralement, on peut
définir des \emph{factorisations par blocs} de $w$ :

\begin{defn}
\label{deffact}
Soit $w\in W$, de longueur $n \geq 1$. Soit $1 \leq p \leq n$. Une
\emph{factorisation en $p$ blocs} de $w$ est un $p$-uplet $(u_1,
\dots, u_p)$ d'éléments de $W - \{1\}$ tels que $w=u_1\dots u_p$, avec
$\ell(w)=\ell(u_1) +\dots + \ell(u_p)$. On note $\Fact_p(w)$ l'ensemble de
ces factorisations.
\end{defn}

Si $w$ est de longueur $n$, on a donc
$\Fact_n(w)=\Red (w)$.

\medskip

Ci-dessous on définit l'action d'Hurwitz du groupe de tresses $B_p$ sur les
factorisations en $p$ blocs d'un élément de $W$. Une factorisation est dite
\emph{primitive} lorsqu'elle ne comporte qu'un seul facteur qui n'est pas
une réflexion (définition \ref{deffactprim}). Dans cet article on détermine les
orbites sous l'action d'Hurwitz des factorisations primitives d'un élément
de Coxeter de $W$ (définition \ref{defcox}).

\subsection*{Action d'Hurwitz}
~ 

Considérons un groupe $G$, et un entier $n\geq 1$. On fait agir le
groupe de tresses à $n$ brins, noté $B_n$, sur le produit cartésien $G^n$ :

\begin{defn}[Action d'Hurwitz]
\label{defhur}
  Soient $\bm{\sigma_1},\dots, \bm{\sigma_{n-1}}$ les générateurs standards de
  $B_n$. L'action d'Hurwitz (à droite) de $B_n$ sur $G^n$ est
  définie par :
  \[ (g_1,\dots, g_n) \cdot \bm{\sigma_i} \ := \ (g_1, \dots, g_{i-1},\
  g_{i+1},\  g_{i+1} ^{-1} g_i g_{i+1},\  g_{i+2}, \dots, g_n) \] pour tout
  $(g_1,\dots, g_n) \in G^n$ et tout $i\in \{1,\dots,n-1\}$.
\end{defn}

On vérifie aisément que cette définition donne bien une action du groupe de
tresses. Cela vient essentiellement du fait que la conjugaison fait partie
de la classe des opérations auto-distributives :
$f*(g*h)=(f*g)*(f*h)$ (ici $f*g=f^{-1}gf$). L'intérêt de ces opérations a
été décrit par Brieskorn (ensembles automorphes \cite{Brieskorn}) et plus
récemment par Dehornoy (LD-systèmes \cite{LD}).

\subsection*{Factorisations primitives}
~ 

L'action d'Hurwitz préserve les fibres de l'application produit $(g_1,
\dots, g_n) \mapsto g_1\dots g_n $. Ainsi, si l'on revient au cadre décrit
plus haut, où $W$ est un groupe de réflexions complexe muni de sa partie
génératrice naturelle $\CR$, on obtient une action d'Hurwitz du groupe de
tresses $B_p$ sur l'ensemble des factorisations en $p$ blocs de $w$
(définition \ref{deffact}), pour tout $w\in W$ et $1\leq p \leq \ell(w)$.

Considérons $(u_1,\dots,u_p) \in \Fact_p(w)$. Le $p$-uplet $(\ell(u_1),\dots,
\ell(u_p))$ définit donc une composition (\ie une partition ordonnée) de
l'entier $n=\ell(w)$. Si $\mu$ est une composition de $n$ en $p$ termes, on note
$\Fact_\mu(w)$ l'ensemble des factorisations dont les longueurs des
facteurs forment la composition $\mu$.

Si l'on fait agir $B_p$ sur un élément de $\Fact_\mu(w)$, la composition de
$n$ associée peut être modifiée, mais pas la partition de $n$ (non
ordonnée) sous-jacente. En effet, l'ensemble $\CR$ des réflexions étant
invariant par conjugaison, la longueur $\ell$ est également invariante par
conjugaison.

On dit qu'une factorisation en $p$ blocs est de forme $\lambda$, pour
$\lambda$ partition de $n$ (notation :
$\lambda \vdash n$), si le $p$-uplet non ordonné des longueurs des
facteurs constitue la partition $\lambda$. On note $\Fact_\lambda (w)$
l'ensemble des factorisations de forme $\lambda$. Ainsi, l'action
d'Hurwitz de $B_p$ sur $\Fact_p (w)$ stabilise les ensembles
$\Fact_\lambda (w)$, pour $\lambda \vdash n$ et $\#\lambda =p$
($\#\lambda$ désignant le nombre de parts de la partition $\lambda$).

Pour $\lambda \vdash n$, on note classiquement $\lambda=k_1^{p_1}\dots
k_r^{p_r}$ (où $k_1 > \dots > k_r \geq1$ et $p_i \geq 1$), si la partition
est constituée de $p_1$ fois l'entier $k_1$, $\dots$, $p_r$ fois l'entier
$k_r$.

\begin{defn}
  \label{deffactprim}
  Une \emph{factorisation primitive} de $w$ est une factorisation par blocs
  de $w$, de forme $k^1 1^{n-k}\vdash n$, pour $k \geq 2$ : tous les blocs
  sauf un doivent être des réflexions. L'unique facteur de longueur
  strictement supérieure à $1$ est appelé \emph{facteur long} de la
  factorisation primitive.
\end{defn}

On étudie ici les factorisations d'un \emph{élément de Coxeter}, au sens de
\cite{BessisKPi1} (voir définition \ref{defcox}). Dans le cas réel, on
retrouve la notion classique d'élément de Coxeter. Dans la suite de
l'introduction, on précisera l'intérêt fondamental des éléments de Coxeter
dans certaines constructions algébriques relatives à un groupe de
réflexions complexe bien engendré.

Considérons un élément de Coxeter $c$ de $W$, et notons $n$ sa
longueur. Le cas le plus élémentaire des factorisations de forme $1^n
\vdash n$, \ie des décompositions réduites de $c$, est connu : on
sait que l'action d'Hurwitz sur $\Red(c)=\Fact_{1^n}(c)$ est
transitive. Dans un premier temps cette propriété a été prouvée pour les
groupes réels (\cite{deligneletter}, puis \cite[Prop.1.6.1]{Bessisdual}),
avec une preuve générale, puis vérifiée au cas par cas pour le reste des
groupes complexes \cite[Prop.7.5]{BessisKPi1}.

Le cas que nous traitons ici peut se voir comme l'étape suivante,
c'est-à-dire la détermination des orbites d'Hurwitz des factorisations
primitives de $c$. Le résultat principal de cet article est le suivant :

\begin{thm}[Orbites d'Hurwitz primitives]
\label{thmintro}
  Soit $W$ un groupe de réflexions complexe fini, bien engendré. Soit $c$
  un élément de Coxeter de $W$. Alors, deux factorisations primitives de
  $c$ sont dans la même orbite d'Hurwitz si et seulement si leurs facteurs
  longs sont conjugués.
\end{thm}

L'idée de la preuve est d'utiliser certaines constructions de Bessis
\cite{BessisKPi1} pour obtenir les factorisations de $c$ à partir de la
géométrie du discriminant de $W$ dans l'espace-quotient $W \qg V$ (qui
s'identifie à $\BC^n$ si on a choisi des invariants fondamentaux $f_1,
\dots, f_n$). On interprète ensuite l'action d'Hurwitz en termes
géométriques.

Dans une première partie on rappelle les définitions classiques et les
propriétés de l'ordre de divisibilité dans le treillis des partitions non
croisées généralisées, noté $\NCP_W(c)$. La partie \ref{partfconj} présente
une notion de forte conjugaison pour les éléments de $\NCP_W(c)$, qui
permet de formuler une version plus forte du théorème principal \ref{thmintro}.

Par la suite, on considère l'hypersurface $\CH$ de $\BC^n$, d'équation le
discriminant de $W$, qui admet une stratification naturelle (par les
orbites de plats). On projette cette hypersurface sur $\BC^{n-1}$ (en
oubliant la coordonnée correspondant à l'invariant de plus haut degré
$f_n$), et on en déduit une stratification du lieu de bifurcation $\CK$ de
$\CH$. En partie \ref{partLL} on rappelle les propriétés du revêtement de
Lyashko-Looijenga $\LL$ , puis dans la partie suivante on modifie légèrement
les constructions de \cite{BessisKPi1}, en définissant une application
$\fact$ qui produit des factorisations associées à la géométrie de
$\CH$. Dans la partie \ref{partLLstrat}, on démontre des propriétés plus
fines du morphisme $\LL$, qui apparaît alors comme un revêtement ramifié
«~stratifié~». On en déduit que l'ensemble des orbites d'Hurwitz des
factorisations primitives est en bijection avec l'ensemble des composantes
connexes par arcs de certains ouverts des strates de $\CK$. La
démonstration du théorème revient ainsi à un problème de connexité dans
$\CK$ (traité en parties \ref{partecp} et \ref{partirred}), et est achevée
en partie \ref{partrefl} par l'étude du cas particulier des réflexions.

Dans la partie \ref{partecp}, on détaille les propriétés de la
stratification de $\CH$. À cette occasion, des conséquences remarquables de
certaines propriétés énoncées dans \cite{BessisKPi1} sont exposées, en
particulier concernant les sous-groupes paraboliques de $W$, leurs classes
de conjugaison et leurs éléments de Coxeter.

\bigskip

La suite de cette introduction présente les motivations du problème.

\subsection*{Treillis des partitions non croisées, monoïde de
  tresses dual}
~ 

Fixons un élément de Coxeter $c$ dans $W$. On note $\<$ l'ordre
partiel sur $W$ associé à la $R$-longueur $\ell$ (voir définition
\ref{defordre}) :
\[u \< v \mathrm{\ si\ et\ seulement\ si} \ell (u) + \ell (u  ^{-1} v) = \ell
(v) .\]

On considère l'ensemble des diviseurs de $c$ : $\NCP_W(c):= \{ w \in W \tq
w\< c \} $ (la notation $\NCP$ signifie «~non-crossing partitions~»,
cf. interprétation ci-dessous).

L'ensemble partiellement ordonné $(\NCP_W (c), \<)$ permet de construire
les générateurs (appelés les \emph{simples}) du \emph{monoïde de tresses
  dual} de $W$ (cf. \cite{BKL} pour le type $A$, \cite{Bessisdual} et
\cite{BWordre2} pour les groupes réels, puis \cite{BessisKPi1} pour le cas
général).

\medskip

L'étude de la combinatoire de $\NCP_W(c)$ est au départ motivée par les
conséquences sur les propriétés algébriques du monoïde dual. Ainsi, une des
propriétés remarquables de $\NCP_W(c)$ est qu'il forme un treillis pour l'ordre
$\<$ : preuve générale dans le cas réel par Brady et Watt \cite{BWtreillis}
(voir aussi \cite{IngThom}), le reste au cas par cas, cf. \cite[Lemme
8.9]{BessisKPi1}. C'est essentiellement cette propriété de treillis qui
donne au monoïde de tresses dual sa structure de monoïde de Garside (comme
défini dans \cite{Deh2}).

Mais la richesse de la structure combinatoire de $\NCP_W(c)$ en a fait plus
récemment un objet d'étude en soi (voir par exemple \cite{Armstrong}). Dans
le cas du groupe $W(A_{n-1})$, ce treillis correspond à l'ensemble des
partitions non croisées d'un $n$-gone, et est à la base de la construction
du monoïde de Birman-Ko-Lee \cite{BKL,BDM}. Pour les groupes de type
$G(e,e,r)$, on peut également interpréter l'ensemble des diviseurs de $c$
comme certains types de partitions non croisées \cite{BessisCorran}. Par
extension, l'ensemble $\NCP_W(c)$ des diviseurs d'un élément de Coxeter
dans un groupe bien engendré est appelé le \emph{treillis des partitions
  non croisées généralisées} (d'où la notation $\NCP$).

\subsection*{\texorpdfstring{Formule de Chapoton pour les chaînes de $\NCP_W(c)$ et
  factorisations de $c$}{Formule de Chapoton pour les chaînes de NCPW(c) et
  factorisations de c}}
~ 

Si $W$ n'est pas irréductible et s'écrit $W_1 \times W_2$, alors le
treillis $\NCP_W(c)$ est le produit des treillis $\NCP_{W_1}(c_1)$ et
$\NCP_{W_2}(c_2)$ (où $c_i$ désigne la composante de $c$ sur $W_i$), muni
de l'ordre produit. Dans cette partie on supposera que $W$ est
irréductible.

Notons $d_1 \leq \dots \leq d_n = h$ les degrés invariants de $W$. Une
formule, énoncée initialement par Chapoton, exprime le nombre de chaînes de
longueur donnée dans le treillis $\NCP_W(c)$ en fonction de ces
degrés. Notons $Z_W$ le polynôme \emph{Zêta} de $\NCP_W(c)$ : c'est la
fonction sur $\BN$ telle que pour tout $N$, $Z_W(N)$ est le nombre de
chaînes larges (ou multi-chaînes) $w_1 \< \dots \< w_{N-1}$ dans
$(\NCP_W(c),\<)$ (de façon générale $Z_W$ est toujours un polynôme : voir
par exemple Stanley \cite[Chap. 3.11]{stanley}).

\begin{prop}
\label{propchapoton}
Soit $W$ un groupe de réflexions complexe bien engendré,
irréductible. Alors, avec les notations ci-dessus, on a :
\[ Z_W (N) = \prod_{i=1}^n \frac{d_i + (N-1)h}{d_i} .\] 
\end{prop}

Dans le cas des groupes réels, cette formule a été observée par Chapoton
\cite[Prop. 9]{chapoton}, utilisant des calculs de Reiner et Athanasiadis
\cite{reiner,athareiner}. La généralisation aux groupes complexes est
énoncée par Bessis \cite[Prop. 13.1]{BessisKPi1} et utilise des résultats
de Bessis et Corran \cite{BessisCorran}.  La preuve est au cas par cas :
méthodes \emph{ad hoc} pour les séries infinies, et logiciel GAP pour les
types exceptionnels. Pour le cas réel, une démonstration plus directe
(utilisant une formule de récurrence générale qui permet de faire la
vérification au cas par cas plus simplement) a été récemment publiée par
Reading \cite{reading}.  L'apparente simplicité de la formule de Chapoton
motive toujours la recherche d'une explication complètement générale.

En posant $N=2$ dans l'expression de $Z_W(N)$ on obtient naturellement le
cardinal de $\NCP_W(c)$ :
\[ \left| \NCP_W (c) \right| = \prod_{i=1}^n \frac{d_i + h}{d_i}. \] Même
pour cette formule, la seule preuve connue à l'heure actuelle est au cas
par cas. L'entier $ \left| \NCP_W (c) \right|$, parfois noté $\Cat(W)$, est
appelé \emph{nombre de Catalan de type $W$} : pour $W=A_{n-1}$ on obtient
en effet le nombre de Catalan classique $\frac{1}{n+1} \binom{2n}{n}$.
 
\medskip

Les factorisations de $c$ sont bien évidemment reliées aux chaînes de
$\NCP_W(c)$. Partant d'une chaîne $w_1 \< ... \< w_p \< c$, on peut obtenir
une factorisation par blocs $c=u_1 \dots u_{p+1}$ en posant $u_1 = w_1$,
$u_i = w_{i-1}^{-1} w_i$ pour $ 2\leq i \leq p$, et $u_{p+1}=w_p ^{-1}
c$. Inversement, à partir d'une telle factorisation de $c$ on peut obtenir
une chaîne de $\NCP_W(c)$. Ainsi, dénombrer les chaînes dans $\NCP_W(c)$
revient à compter le nombre de \emph{factorisations par blocs} de $c$ de
tailles fixées.

Toutefois, on a choisi ici de travailler sur des factorisations
«~strictes~» (pas de facteurs triviaux), alors que la formule de Chapoton
concerne les chaînes larges. Ce choix est motivé par le fait que ce sont
des factorisations strictes qui apparaissent géométriquement via
l'application «~$\fact$~» (construite en partie \ref{partlbl}). Les nombres de
factorisations strictes et les nombres de chaînes larges sont alors liés
par la formule suivante (on utilise les formules de passage entre nombres
de chaînes strictes et nombres de chaînes larges, données par exemple dans
\cite[Chap. 3.11]{stanley}) :

\[ Z_W (N) = \sum_{k\geq 1} \left| \Fact_k(c) \right| \binom{N}{k} . \]

Une piste pour une démonstration générale de la formule de Chapoton est donc
l'étude de la combinatoire des factorisations par blocs.

\subsection*{\texorpdfstring{Remerciements.}{Remerciements}}
Je tiens à remercier vivement David Bessis pour ses conseils et son aide
sur de nombreux points de cet article. Merci également à Emmanuel Lepage
dont une remarque judicieuse a permis de simplifier la partie
\ref{subpartcpa}.

\section{\texorpdfstring{Ordre de divisibilité dans $\NCP_W$}{Ordre de divisibilité dans NCPW}}
\label{partordre}

\subsection{Éléments de Coxeter}
\label{subpartcox}
~ 

Soit $V$ un espace vectoriel complexe de dimension $n\geq 1$, et $W
\subseteq \GL(V)$ un groupe de réflexions complexe. On suppose que $W$ est
essentiel, \ie que $V^W=\{0\}$. Dans ce cas, on dit que $W$ est un groupe
de réflexions \emph{bien engendré} si un ensemble de $n$ réflexions suffit
à l'engendrer. Désormais on supposera toujours que $W$ est un groupe de
réflexions complexe bien engendré.

Soit $\alpha$ une racine de l'unité. On dit qu'un élément de $W$ est
\emph{$\alpha$-régulier} (au sens de Springer \cite{Springer}) s'il possède
un vecteur propre régulier (\ie n'appartenant à aucun des hyperplans de
réflexions) pour la valeur propre $\alpha$. Si $W$ est irréductible, notons
$d_1 \leq \dots \leq d_n$ les degrés invariants de $W$. On note $h$ et on
appelle \emph{nombre de Coxeter} le plus grand degré $d_n$. D'après la
théorie de Springer (\cite{Springer}, résumé dans
\cite[Thm. 1.9]{BessisKPi1}), il existe dans $W$ des éléments
$\eh$-réguliers. Dans ce cadre, en généralisant le cas réel, on pose la
définition suivante.

\begin{defn}
  \label{defcox}
  Soit $W$ un groupe de réflexions complexe fini bien engendré. Si $W$ est
  irréductible, un \emph{élément de Coxeter} de $W$ est un élément
  $\eh$-régulier. Dans le cas général, on appelle
  \emph{élément de Coxeter} un produit d'éléments de Coxeter des
  composantes irréductibles de $W$.
\end{defn} 

\subsection{Ordre de divisibilité dans $W$}
\label{subpartordre}
~ 

On définit ci-dessous l'ordre de divisibilité à gauche dans $W$ (parfois
appelé ordre absolu). Rappelons que pour $w \in W$, $\ell(w)$ désigne la
longueur d'une décomposition réduite de $w$ en réflexions de $W$.

\begin{defn}
  \label{defordre}
  Soient $u,v \in W$. On pose :
  \[u \< v \mathrm{\ si\ et\ seulement\ si} \ell (u) + \ell (u  ^{-1} v) = \ell
(v)\]
\end{defn}

La relation $\<$ munit $W$ d'une structure d'ensemble partiellement ordonné
: $u$ divise~$v$ si et seulement si $u$ peut s'écrire comme un préfixe
d'une décomposition réduite de $v$. Comme $\CR$ est invariant par
conjugaison, on a aussi $u \< v$ si et seulement si $u$ est un suffixe (ou
même un sous-facteur) d'une décomposition réduite de $v$. Ainsi, l'ordre
$\<$ coïncide avec l'ordre de divisibilité à droite dans $W$.

\medskip

Considérons un élément de Coxeter $c$ de $W$. On note $\NCP_W(c)$
l'ensemble des diviseurs de $c$ dans $W$. D'après la théorie de Springer,
tous les éléments de Coxeter de $W$ sont conjugués. D'autre part, l'ordre
$\<$ est invariant par conjugaison, donc si $c'=aca^{-1}$, alors
$\NCP_W(c')=a\NCP_W(c)a^{-1}$. Ainsi, la structure de $\NCP_W(c)$ ne dépend
pas du choix de l'élément de Coxeter $c$ ; on notera souvent simplement
$\NCP_W$.

\medskip

Dans $\NCP_W$, la longueur correspond toujours à la codimension de l'espace
des points fixes :

\begin{prop}
  \label{proplongueur}
  Soit $w\in \NCP_W$. Alors :
  \[ \ell(w) = \codim \Ker(w-1) .\]
\end{prop}

\begin{rqe}
  Cela revient à dire que pour les éléments de $\NCP_W$, la longueur
  relativement à $\CR$ est la même que la longueur relativement à
  l'ensemble de toutes les réflexions du groupe unitaire $\U(V)$. Dans le
  cas des groupes réels, cette propriété est vraie pour tout élément de $W$
  (cf. \cite[Lemme 2.8]{carter}, ou \cite[Prop.2.2]{BWordre2}).
\end{rqe}

\begin{proof}[Démonstration :]
  On a toujours $\ell(w) \geq \codim \Ker(w-1)$. En effet, si $w=r_1\dots
  r_k$ alors $\Ker(w-1) \supseteq {\bigcap_i \Ker(r_i-1)}$. Donc pour
  il suffit de montrer que $\ell(c)=n$ (considérer $w$ et
  $w^{-1}c$). L'inégalité $\ell(c) \geq n$ est évidente car $\Ker(c-1)=\{0\}$
  d'après la théorie de Springer. On conclut en notant qu'on peut
  construire géométriquement des factorisations de $c$ en $n$ réflexions
  (\cite[Lemme 7.2]{BessisKPi1} ou partie \ref{partlbl} de cet article).
\end{proof}

\subsection{Théorème de Brady-Watt} 
\label{subpartbw}
~ 

Une conséquence de la proposition \ref{proplongueur} est que l'ordre~$\<$
dans $\NCP_W$ est la restriction de l'ordre partiel sur $\U(V)$ étudié par
Brady et Watt dans \cite{BWordre} : $u \< v$ si et seulement si $\Im(u-1)
\oplus \Im(u^{-1}v -1)= \Im(v)$. On en déduit en particulier :

\begin{prop}
  \label{propordre}
  Soit $w \in \NCP_W$. Si $(r_1,\dots, r_p)$ est une décomposition réduite
  de $w$, alors $\Ker(w-1)= \Ker(r_1-1) \cap \dots \cap \Ker(r_k-1)$.
\end{prop}

Les résultats de \cite{BWordre} sont donnés dans le cadre de $\OO(V)$ pour
$V$ espace euclidien, mais restent valables sans modifications dans $\U(V)$
avec $V$ hermitien complexe. On pourra donc appliquer le théorème suivant :

\begin{thm}[Brady-Watt]
  \label{thmBW}
  Pour tout $g \in \U(V)$, l'application
  \[ \begin{array}{rcl} \left(\{ f \in \U(V) \tq f \< g \},\ \< \right)
    &\rightarrow & \left(\{\mathrm{s.e.v.\ de\ } V \mathrm{\ contenant\ }
      \Ker(g-1) \},\ \supseteq \right)\\ f & \mapsto & \Ker(f-1)
    \end{array} \]
    est un isomorphisme d'ensembles partiellement ordonnés.
  \end{thm}

Par conséquent, l'application $w \mapsto \Ker(w-1)$ est injective sur $
\NCP_W $.

\section{\texorpdfstring{Action d'Hurwitz sur les factorisations primitives
    et conjugaison forte dans $\NCP_W$}{Action d'Hurwitz sur les
    factorisations primitives et conjugaison forte}}
\label{partfconj}

Soit $p \in \{1,\dots, n\}$. Le groupe de tresses classique $B_p$ agit par
action d'Hurwitz sur l'ensemble $\Fact _p (c)$ des factorisations en $p$
blocs de $c$, selon la définition \ref{defhur}.
Dans le cas où $p=n$, l'action de $B_n$ est transitive sur $\Fact_n (c)
=\Red(c)$. Ce n'est pas vrai dans le cas général, puisque des invariants
évidents apparaissent :
\begin{itemize}
\item la partition de $n$ associée au $p$-uplet des longueurs des facteurs ;
\item le multi-ensemble des classes de conjugaison des facteurs (car
  l'action d'Hurwitz conjugue les facteurs).
\end{itemize}

Plus précisément, l'action d'Hurwitz conjugue les facteurs par d'autres
facteurs «~dont les longueurs s'ajoutent~». Cela motive l'introduction
d'une notion plus forte de conjugaison.

\begin{defn}
\label{deffconj}
Pour $w,w' \in \NCP_W(c)$, on définit la relation «~$w$ et $w'$ sont
fortement conjugués dans $\NCP_W(c)$~» (notée $w \fconj w'$) comme la
clôture transitive et symétrique de la relation suivante :
\[w \stackrel{1}{\sim} w' \text{\ s'il\ existe\ } x \in W \text{\ tel\ que\
} xw=w'x , \text{\ avec\ } xw \in \NCP_W(c) \text{\ et\ } \ell(xw)=\ell(x)
+ \ell(w) .\]
\end{defn}

\begin{rqe}
\label{rqhur}
  Bien sûr la relation $ \stackrel{1}{\sim}$ reflète l'incidence, sur un
  facteur, de l'action d'Hurwitz d'une tresse élémentaire $\bm{\sigma_i}$
  sur une factorisation de $c$. Par transitivité, si une tresse $\beta$
  transforme une factorisation $\xi$ en une factorisation $\xi'$, et envoie
  le numéro d'un facteur $w$ de $\xi$ sur celui d'un facteur $w'$ de
  $\xi'$, alors $w$ et $w'$ sont fortement conjugués.
\end{rqe}

On voit facilement qu'on obtient la même relation d'équivalence si on
impose que les conjugateurs $x$ de la définition soient toujours des
réflexions de $\NCP_W$. En effet, si $w$ et $w'$ sont tels que $w
\stackrel{1}{\sim} w'$, notons $x$ un élément de $W$ tel que $xw=w'x
\in \NCP_W$, avec $\ell(xw)=\ell(x) + \ell(w)$. Si l'on décompose $x$ en
produit minimal de réflexions ($x=r_1\dots r_k$), et si l'on note $u
\stackrel{\mathrm{ref}}{\sim} v$ la relation «~$u$ est élémentairement
conjugué à $v$ par une réflexion~», alors on a : $w
\stackrel{\mathrm{ref}}{\sim} r_k w r_k ^{-1} \stackrel{\mathrm{ref}}{\sim}
r_{k-1} r_k w r_k ^{-1} r_{k-1}^{-1} \stackrel{\mathrm{ref}}{\sim} \dots
\stackrel{\mathrm{ref}}{\sim} xwx^{-1}=w'$.
Cette remarque permet de faire le lien avec les orbites d'Hurwitz
de \emph{factorisations primitives} de $c$, où tous les facteurs sauf un
sont des réflexions (cf. définition \ref{deffactprim}).

\begin{prop}
  \label{propfconjhur}
  Soient $u,v \in \NCP_W$, de longueurs strictement supérieures à
  $1$. Alors les propriétés suivantes sont équivalentes :
  \begin{enumerate}[(i)]
  \item $u$ et $v$ sont fortement conjugués dans $\NCP_W$ ;
  \item il existe $\xi=(u,r_{2},\dots,r_k)$ et
    $\xi'=(v,r_{2}',\dots,r_k')$ deux factorisations primitives de
    $c$, avec $u$ (resp. $v$) facteur long de $\xi$ (resp. $\xi'$),
    telles que $\xi$ et $\xi'$ soient dans la même orbite d'Hurwitz
    sous $B_k$.
  \end{enumerate}
\end{prop}

\begin{rqe} \label{rqrefl} Dans le cas des réflexions, on a une proposition
  similaire, mais il faut s'assurer que la (permutation associée à la)
  tresse qui transforme $\xi$ en $\xi'$ envoie le numéro de la
  position de $u$ sur celui de la position de $v$.
\end{rqe}

\begin{proof}[Démonstration :]
  (ii) $\Rightarrow$ (i) est clair par définition, puisque, étant donnée
  l'invariance de la longueur par conjugaison, $v$ est nécessairement
  obtenu à partir de $u$ par forte conjugaison (cf. remarque
  \ref{rqhur}). Pour le sens direct, montrons d'abord (ii) lorsque $u$ et
  $v$ sont tels que $ur=rv \< c$, avec $r \in \CR$ et
  $\ell(ur)=\ell(u)+1$. Il suffit dans ce cas de prendre une factorisation
  quelconque de $c$ qui commence par $(u,r)$, et de faire agir la tresse
  $\bm{\sigma_1} ^2$. On conclut par transitivité.
\end{proof}

Ainsi le théorème \ref{thmintro} se déduit du théorème suivant, qui donne
également des propriétés supplémentaires de l'action d'Hurwitz sur
$\Red(c)$ :

\begin{thm}
  \label{thmfconj}
  Soit $u,v \in \NCP_W$. Alors $u$ et $v$ sont fortement conjugués dans
  $\NCP_W$ si et seulement si $u$ et $v$ sont conjugués dans $W$.
\end{thm}

On aura ainsi déterminé les orbites d'Hurwitz de factorisations
primitives. Toutefois, cela ne suffit pas pour comprendre complètement
l'action d'Hurwitz sur les factorisations de forme quelconque. En effet,
pour toute factorisation, le multi-ensemble des classes de conjugaison
forte des facteurs est naturellement invariant par l'action
d'Hurwitz. Cependant, pour les factorisations non primitives, la condition
n'est en général pas suffisante pour que deux factorisations soient dans la
même orbite. Par exemple, dans le cas où $p=2$, l'orbite d'Hurwitz de
$(u_1,u_2) \in \Fact_2 (c)$ est $\{ (u_1^{c^k},u_2^{c^k}),
(u_2^{c^{k+1}},u_1^{c^k}), \linebreak[1] {k \in \BZ} \}$ (où la notation
$u^v$ désigne le conjugué $v^{-1}uv$). Donc dans le type $A$, cela revient
à agir par rotation sur le diagramme des partitions non-croisées. On peut
ainsi facilement trouver un contre-exemple dans $\FS_6$. Posons $u_1=(2\
3)(1\ 5\ 6), u_2=(1\ 3\ 4)$, et $v_1=(3\ 4)(1\ 5\ 6)$, $v_2=(1\ 2\ 4)$ :
alors les factorisations $(u_1,u_2)$ et $(v_1,v_2)$ de $\Fact_2((1\ 2\ 3\
4\ 5\ 6))$ ne sont pas dans la même orbite d'Hurwitz sous $B_2$, alors que
$u_1$ et $v_1$, resp. $u_2$ et $v_2$, sont conjugués. Le problème vient
du fait que lorsqu'on fait agir $B_p$ par action d'Hurwitz, les
conjugateurs qui s'appliquent ne sont pas quelconques, mais doivent être
des éléments qui sont aussi des facteurs (en particulier pas nécessairement
des réflexions).

\medskip

Désormais on supposera toujours que $W$ est irréductible (et bien
engendré). Si $W$ n'est pas irréductible, le théorème \ref{thmfconj} pour
$W$ se déduit du cas irréductible : on vérifie aisément qu'il suffit
d'appliquer le résultat à chacune des composantes irréductibles de $W$.

\section{Le morphisme de Lyashko-Looijenga}
\label{partLL}
Dans cette partie on rappelle la construction du morphisme de
Lyashko-Looijenga de \cite[Part. 5]{BessisKPi1}, et on fixe des notations
et conventions concernant l'espace de configuration de $n$ points.

\medskip

Le morphisme de Lyashko-Looijenga a été introduit par Lyashko en 1973
(selon Arnold \cite{arnold}) et indépendamment par Looijenga dans
\cite{looijenga} : voir \cite[Chap.5.1]{LZgraphs} et \cite{landozvonkine}
pour un historique détaillé.  Bessis, dans \cite{BessisKPi1}, a généralisé la
définition de $\LL$ à tous les groupes de réflexions complexes bien
engendrés, le cas inititial correspondant aux groupes de Weyl.

\subsection{Discriminant d'un groupe bien engendré}
\label{subpartLL0}
~

Soit $W \subseteq \GL(V)$ un groupe de réflexions complexe bien engendré,
irréductible. Fixons une base $v_1,\dots, v_n$ de $V^*$, et $f_1,\dots,
f_n$ dans $S(V^*)= \BC[v_1,\dots, v_n]$ un système d'invariants
fondamentaux, polynômes homogènes de degrés uniquement déterminés $d_1 \leq
\dots \leq d_n=h$, tels que $\BC[v_1,\dots, v_n]^W=\BC[f_1,\dots, f_n]$
(théorème de Shephard-Todd-Chevalley). On a l'isomorphisme
\[\begin{array}{rcl}
W \qg V & \overset\sim\to & \BC^n \\
\bar{v} & \mapsto & (f_1(v),\dots, f_n(v)).
\end{array}
\]

Notons $\CA$ l'ensemble des hyperplans de
 réflexions de $W$.  L'équation de $\bigcup_{H \in \CA} H$ peut s'écrire
 comme un polynôme invariant, noté $\Delta$, le discriminant de $W$ :

\[ \Delta = \prod_{H\in \CA} \alpha_H ^{e_H} \in \BC[f_1,\dots,
f_n], \]
où $\alpha_H$ est une forme linéaire de noyau $H$ et $e_H$ est l'ordre du
sous-groupe parabolique cyclique $W_H$.

Lorsque $W$ est bien engendré, il existe un système d'invariants
$f_1,\dots, f_n$ tel que le discriminant $\Delta$ de $W$ s'écrive :
\[\Delta= f_n^n + a_2 f_n ^{n-2} +\dots + a_n ,\]
où les $a_i$ sont des polynômes en $f_1, \dots, f_{n-1}$
\cite[Thm. 2.4]{BessisKPi1}. La propriété fondamentale est que $\Delta$ est
un polynôme monique de degré $n$ en $f_n$ (le coefficient de $f_n^{n-1}$
est rendu nul par simple translation de la variable $f_n$).

On note :
\[\CH := \{\bar{v}\in W \qg V \tq \D(\bar{v})=0 \} = p\left(\bigcup_{H \in \CA}
H\right) ,\] où $p$ est le morphisme quotient $V \surj \ W \qg V$. Désormais,
lorsqu'on se placera dans l'espace quotient, on considèrera $(f_1,\dots,
f_n)$ comme les coefficients des points de $W \qg V \simeq \BC^n$.

Le morphisme de Lyashko-Looijenga $\LL$, que l'on définira
précisément en partie \ref{subpartLL}, permet d'étudier l'hypersurface $\CH$ à
travers les fibres de la projection $(f_1,\dots,f_n) \mapsto (f_1,\dots,
f_{n-1})$. Ensemblistement, $\LL$ associe à $(f_1,\dots, f_{n-1}) \in
\BC^{n-1}$ le multi-ensemble des racines de $\D (f_1,\dots, f_n)$
vu comme polynôme en $f_n$. Dans la partie suivante on fixe les notations
et définitions concernant l'espace d'arrivée de $\LL$.

\subsection{L'espace des configurations de $n$ points}
\label{subpartEn}
~

On note $E_n$ l'ensemble des configurations
centrées de $n$ points, avec multiplicités, \ie
\[E_n := H_0/ \FS_n\] où $H_0$ est l'hyperplan de $\BC^n$ d'équation $\sum
x_i =0$. Considérons le morphisme quotient ${f: H_0 \to E_n}$. C'est un
morphisme fini (donc fermé), correspondant à l'inclusion
\[\BC[e_1,\dots,e_n]/(e_1)\ \inj \
\BC[x_1,\dots,x_n]/(\sum x_i),\] où $e_1,\dots,e_n$ désignent les fonctions
symétriques élémentaires en les $x_i$. L'espace $E_n$ est une variété
algébrique, son anneau des fonctions régulières sera noté
$\BC[e_2,\dots,e_n]$.

\medskip

On peut stratifier naturellement $E_n$ par les partitions de l'entier
$n$. Pour $\lambda \vdash n$, on note $E_{\lambda}^0$ l'ensemble des
configurations $X\in E_n$ dont les multiplicités sont distribuées selon
$\lambda$. Ainsi, $E_n$ est l'union disjointe des $E_{\lambda}^0$.

Une partie $E_{\lambda}^0$ n'est pas fermée ;
des points peuvent fusionner et on obtiendra une configuration
correspondant à une partition moins fine :

\begin{defn}
  Soient $\lambda, \mu$ deux partitions de $n$. On note $\mu \leq \lambda$
  la relation «~$\mu$ est \emph{moins fine} que $\lambda$~», \ie $\mu$ est
  obtenue à partir de $\lambda$ après un ou plusieurs regroupements de
  parts.

  On pose $E_{\lambda} :=\bigsqcup_{\mu \leq \lambda} E_{\mu}^0$
  (le symbole $\bigsqcup$ désigne une union disjointe).
\end{defn}

Ainsi, pour $\eps:=1^n \vdash n$, $E_{\eps}$ est l'espace $E_n$
entier, et $E_{\eps}^0$ est l'ensemble des points réguliers de $E_n$,
que l'on notera $\Enreg$. On a $\Enreg=H_0^{\mathrm{reg}}/\FS_n$, où
$H_0^{\mathrm{reg}}=\{(x_1,\dots,x_n)\in H_0 \tq \forall i\neq j,\ x_i\neq
x_j \}$. Pour $\alpha:=2^1 1^{n-2} \vdash n$, on a $E_{\alpha}=E_n-\Enreg$.

Pour la topologie classique, comme pour la topologie de Zariski,
$E_{\lambda}$ est un fermé qui est l'adhérence de
$E_{\lambda}^0$.

Fixons une configuration de référence $X^\circ$ dans $\Enreg$, par exemple
sur la droite réelle. On note $B_n:=\pi_1(\Enreg,X^\circ)$ le groupe de
tresses à $n$ brins. On a la présentation classique :
\[B_n  \simeq \left< \bm{\sigma_1},\dots,\bm{\sigma_{n-1}}
 \left| \bm{\sigma_i}\bm{\sigma_{i+1}} \bm{\sigma_i} =
\bm{\sigma_{i+1}}\bm{\sigma_i} \bm{\sigma_{i+1}},
\bm{\sigma_i}\bm{\sigma_j} = \bm{\sigma_j}\bm{\sigma_i} \text{ pour } |i-j| >
1 \right. \right> ,\] 
où par convention $\bm{\sigma_i}$ est représenté par le chemin suivant dans
le plan complexe :

{\shorthandoff{?;:}%
\[\xy
(-10,0)="1", (-2,0)="2", (7,0)="3", 
(14,0)="4", (22,0)="5", (32,0)="6", (40,0)="7",
"1"*{\bullet},"2"*{\bullet},"3"*{\bullet},"4"*{\bullet},
"5"*{\bullet},"6"*{\bullet},"7"*{\bullet},
(14,-3)*{_{x_i}},(24,-3)*{_{x_{i+1}}},
(-10,-3)*{_{x_1}},(40,-3)*{_{x_{n}}},
"4";"5" **\crv{(18,-3)}
 ?(.5)*\dir{>},
"5";"4" **\crv{(18,3)}
 ?(.5)*\dir{>}
\endxy \]
}

Tout chemin dans $\Enreg$ peut être vu comme un élément du groupe de
tresses $B_n$. Précisons la construction, qui sera importante par la
suite. On définit l'ordre lexicographique sur $\BC$ : $z \lex z'\ \ssi \
[ \re(z) \leq \re(z')] \mathrm{\ et\ } [\re(z)=\re(z')\ \Rightarrow \ \im(z)
\leq \im(z')] $.

 Soit $X$ une configuration quelconque de $\Enreg$, et posons
$\widetilde{X}=(x_1,\dots,x_n )$ son support ordonné pour
$\lex$. Considérons un chemin $t\mapsto (x_1(t),\dots, x_n(t))$ dans
$H_0^{\mathrm{reg}}$ de $\widetilde{X}$ vers $\widetilde{X^\circ}$ (le support
ordonné de $X^\circ$), tel que pour tout $t$, $x_1(t)\lex \dots \lex
x_n(t)$. Il détermine une unique classe d'homotopie de chemin de $X$ vers
$X^\circ$, que l'on note $\gamma_X$.  Si $X,X' \in \Enreg$, on fixe ainsi
une bijection entre le groupe $\pi_1(\Enreg, X^\circ)$ et l'ensemble des
classes d'homotopie de chemin de $X$ vers $X'$ par : $\tau \mapsto
\gamma_X \cdot \tau \cdot \gamma_{X'}^{-1}$.

Via ces conventions, on pourra désormais considérer tout chemin dans
$\Enreg$ comme un élément de $B_n$. De la même façon, si $\lambda$ est une
partition de $n$, on peut associer à tout chemin de $E_\lambda ^0$ un
élément de $B_r$, où $r=\#\lambda$ est le nombre de parts de $\lambda$ : on
considère simplement le mouvement du support de la configuration.  On
utilisera ces conventions dans la proposition \ref{propcompat} et le lemme
\ref{lemcompat}.

\subsection{\texorpdfstring{Définitions et propriétés du morphisme
    $\LL$}{Définitions et propriétés du morphisme LL}}
\label{subpartLL}
~ 

On pose $Y:=\Spec \BC[f_1,\dots, f_{n-1} ] \simeq \BC^{n-1}$. D'autre part,
$E_n=\Spec \BC[e_2,\dots, e_n]$. Rappelons que $a_2,\dots,a_n$, éléments de
$\BC[f_1,\dots,f_{n-1}]$, désignent les coefficients du discriminant
$\Delta$ en tant que polynôme en $f_n$.

\begin{defn}[{\cite[Def.5.1]{BessisKPi1}}]

  \label{defLL}
  Le morphisme de Lyashko-Looijenga (généralisé) est le morphisme de $Y$
  dans $E_n$ défini algébriquement par :
  \[ \begin{array}{rcl}
    \BC[e_2,\dots,e_n] & \to & \BC[f_1,\dots,f_{n-1}]\\
    e_i & \mapsto & (-1)^i a_i
  \end{array}\]
\end{defn}

Ensemblistement, $\LL$ envoie $y=(f_1,\dots,f_{n-1}) \in Y$ sur le
multi-ensemble (élément de $E_n$) des racines du polynôme $\D = f_n^n + a_2
(y) f_n^{n-2} +\dots + a_n(y)$, \ie les intersections, avec multiplicité,
de $\CH$ avec la droite $L_y:=\{(y,f_n) \tq f_n \in \BC \}$. On peut voir
aussi $\LL$ comme le morphisme algébrique qui envoie $y \in Y$ sur $(a_2,
\dots, a_{n})$ : il est alors quasi-homogène pour les poids : $\deg f_i =
d_i$, $\deg a_i = ih$.

\medskip

On pose $\CK := \{ y\in Y \tq \Disc (\D(y,f_n); f_n) =0 \}$, appelé lieu de
bifurcation de $\D$.  C'est le lieu où $\D$ a des racines multiples en tant
que polynôme en $f_n$. Ainsi $\CK=\LL^{-1}(E_\alpha)$ et $Y-\CK= \LL ^{-1}
(\Enreg)$, avec les notations de la partie \ref{subpartEn}. Les propriétés
fondamentales de $\LL$ sont données par le théorème suivant :

\begin{thm}[{\cite[Thm. 5.3]{BessisKPi1}}]
  \label{thmLLBessis}
  Les polynômes $a_2,\dots,a_n \in \BC[f_1,\dots,f_{n-1}]$ sont
  algébriquement indépendants, et $\BC[f_1,\dots,f_{n-1}]$ est un
  $\BC[a_2,\dots,a_n]$-module libre gradué de rang $\frac{n!h^n}{|W|}$. Par
  conséquent, $LL$ est un morphisme fini. De plus, sa restriction $Y-\CK \surj
  \Enreg$ est un revêtement non ramifié de degré $\frac{n!h^n}{|W|}$.
\end{thm}

La dernière propriété permet de définir, pour chaque $y\in Y - \CK$, une
action de Galois (ou action de monodromie) de $\pi_1(\Enreg,\LL(y))$ sur la
fibre de $\LL(y)$. D'après les conventions de la partie \ref{subpartEn},
cela détermine une action du groupe de tresses $B_n$ sur l'espace $Y- \CK$ :
les orbites sont exactement les fibres de $\LL$ (car $Y-\CK$ est connexe
par arcs). Pour $\beta \in B_n$, on note $y\cdot \beta$ l'image de $y$ par
l'action de $\beta$.
\medskip

Dans la partie suivante, on va construire pour chaque $y\in Y$ une
factorisation de $c$ ; on verra plus loin que l'action d'Hurwitz sur les
factorisations (définition \ref{defhur}) est compatible avec l'action de Galois
définie ici.

\section{Factorisations géométriques}

\label{partlbl}

En adaptant \cite[Part.6]{BessisKPi1}, on construit pour chaque $y\in Y$
une factorisation par blocs de $c$, dont la partition associée correspond à
la distribution des multiplicités de la configuration $\LL(y)$.

\subsection{Tunnels}
~ 

Soit $y \in Y$. Rappelons que $L_y$ désigne la droite complexe $\{(y,x) \tq
x \in \BC \}$. Notons $U_y$ le complémentaire dans $L_y$ des demi-droites
verticales situées au-dessous des points de $\LL(y)$, \ie :
\[ U_y:= \{ (y,z) \in L_y \tq \forall x \in \LL(y),\ \re(z)= \re(x)
\Rightarrow \im(z) > \im (x) \} \]

La partie $\CU := \bigcup_{y\in Y} U_y$ est un ouvert dense et contractile
de $W \qg \Vreg$ \cite[Lemme 6.2]{BessisKPi1}. On peut donc utiliser $\CU$
comme «~point-base~» pour le groupe fondamental de $W \qg \Vreg$. On
définit ainsi le groupe de tresses de $W$ :
\[B(W):=\pi_1(W \qg \Vreg, \CU) \]

\medskip

On rappelle ci-dessous \cite[Def.6.5]{BessisKPi1} :

\begin{defn}
  Un \emph{tunnel} est un triplet $T=(y,z,L)$ tel que $(z,y)\in \CU$,
  $(z+L,y)\in \CU$, et le segment $[(z,y),(z+L,y)]$ est contenu dans $W \qg
  \Vreg$. On associe à $T$ l'élément $b_T$ du groupe de tresses $\pi_1(W
  \qg \Vreg, \CU)$, représenté par le chemin horizontal $t\mapsto
  (z+tL,y)$.
\end{defn}

{\shorthandoff{;:}%
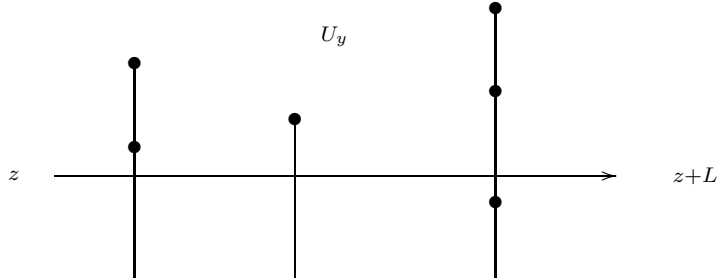
 \[\xy
<5pt,0pt>:<0pt,3.5pt>::
(0,0)="A", (0,-17)="A1", 
(-12,6)="B", (-12,-3)="B1", (-12,-17)="B2", 
(15,12)="C", (15,3)="C1", (15,-9)="C2", (15,-17)="C3", 
(3,9)*{_{U_y}},
"B"*{\bullet},"B1"*{\bullet},"A"*{\bullet},"C"*{\bullet},"C1"*{\bullet},"C2"*{\bullet},
"B";"B1" **@{-},
"B2";"B1" **@{-},
"A";"A1" **@{-},
"C";"C1" **@{-},
"C1";"C2" **@{-},
"C2";"C3" **@{-},
(-18,-6);(24,-6) **@{-}, *\dir{>},
(-21,-6)*{_{z}}, (30,-6)*{_{z+L}}
\endxy \]
}

Soit $y\in Y$, et $(x_1,\dots,x_k)$ le support du multi-ensemble $\LL(y)$,
ordonné selon l'ordre lexicographique $\lex$. Notons $\pr_{\BC}$ la
projection de $W \qg V \simeq Y \times \BC$ sur la dernière
coordonnée. Considérons l'espace $\pr_{\BC} (L_y - U_y) - \{x_1,\dots, x_k
\}$. C'est une union disjointe d'intervalles ouverts, bornés ou non. Pour
$j \in \{1,\dots,k\}$, on note $I_j$ l'intervalle situé sous $x_j$.

{\shorthandoff{;:}%
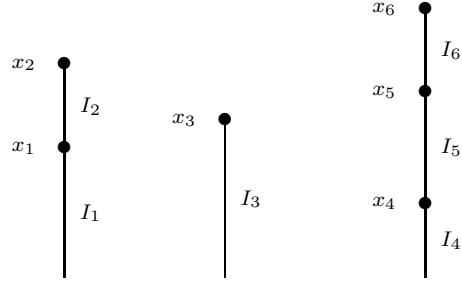
 \[\xy
<5pt,0pt>:<0pt,3.5pt>::
(0,0)="A", (0,-17)="A1", 
(-12,6)="B", (-12,-3)="B1", (-12,-17)="B2", 
(15,12)="C", (15,3)="C1", (15,-9)="C2", (15,-17)="C3", 
"B"*{\bullet},"B1"*{\bullet},"A"*{\bullet},"C"*{\bullet},"C1"*{\bullet},"C2"*{\bullet},
(-15,6)*{_{x_2}},(-15,-3)*{_{x_1}},(-3,0)*{_{x_3}},(12,-9)*{_{x_4}},(12,3)*{_{x_5}},(12,12)*{_{x_6}},
(-10,-10)*{_{I_1}},(-10,1.5)*{_{I_2}},(2,-8.5)*{_{I_3}},(17,-13)*{_{I_4}},(17,-3)*{_{I_5}},(17,7.5)*{_{I_6}},
"B";"B1" **@{-},
"B2";"B1" **@{-},
"A";"A1" **@{-},
"C";"C1" **@{-},
"C1";"C2" **@{-},
"C2";"C3" **@{-}
\endxy \]
}
Pour chaque $x_j$, on choisit un tunnel $T_j$ qui traverse
$I_j$ et pas les autres intervalles. On note $s_j:=b_{T_j}$ l'élément
de $B(W)$ associé.
{\shorthandoff{;:}%
 \[\xy
<5pt,0pt>:<0pt,3.5pt>::
(0,0)="A", (0,-17)="A1", 
(-12,6)="B", (-12,-3)="B1", (-12,-17)="B2", 
(15,12)="C", (15,3)="C1", (15,-9)="C2", (15,-17)="C3", 
"B"*{\bullet},"B1"*{\bullet},"A"*{\bullet},"C"*{\bullet},"C1"*{\bullet},"C2"*{\bullet},
(-15,6)*{_{x_2}},(-15,-3)*{_{x_1}},(-3,0)*{_{x_3}},(12,-9)*{_{x_4}},(12,3)*{_{x_5}},(12,12)*{_{x_6}},
(-8,-12)*{_{s_1}},(-8,-0.5)*{_{s_2}},(4,-10.5)*{_{s_3}},(19,-15)*{_{s_4}},(19,-5)*{_{s_5}},(19,5.5)*{_{s_6}},
(-16,-10);(-8,-10) **@{-}, *\dir{>},(-16,1.5);(-8,1.5) **@{-}, *\dir{>},
(-4,-8.5);(4,-8.5) **@{-}, *\dir{>},
(11,-13);(19,-13) **@{-}, *\dir{>},(11,-3);(19,-3) **@{-},
*\dir{>},(11,7.5);(19,7.5) **@{-}, *\dir{>},
"B";"B1" **@{-},
"B2";"B1" **@{-},
"A";"A1" **@{-},
"C";"C1" **@{-},
"C1";"C2" **@{-},
"C2";"C3" **@{-}
\endxy \]
}

Dans \cite[Def.6.7]{BessisKPi1}, le $k$-uplet $(s_1,\dots,s_k)$ est
 appelé \emph{label} de y et noté $\lbl(y)$. Ici on va associer à $y$ un
 $k$-uplet légèrement différent, mieux adapté au problème. 

\subsection{Factorisations géométriques}
\label{subpartfactogeom}

\begin{defn}
  \label{deflbl}
  Soit $y\in Y$, $x_1,\dots,x_k$ et $s_1,\dots,s_k$ définis comme
  ci-dessus. On note $\fact_B(y)$ le $k$-uplet $(s_1',\dots, s_k')$
  d'éléments de $B(W)$ où :
  \[ s_i'=\left\{
    \begin{array}{ll}
      s_i s_{i+1}^{-1} & \mathrm{si\ }\re(x_{i+1})=\re(x_i)\\
      s_i & \mathrm{sinon.}
    \end{array}
  \right.\]

  On appelle \emph{factorisation associée à $y$}, et on note $\fact(y)$, le
  $k$-uplet $(\pi(s_1'),\dots, \pi(s_k'))$ d'éléments de $W$ où $\pi:B(W)
  \surj W$.
\end{defn}

La raison de l'utilisation du terme «~factorisation~» sera clarifiée plus
loin.

\begin{rqe}
  \label{rqpi}
  L'application $\pi$ est une projection naturelle $B(W)
  \stackrel{\pi}{\surj}W$. Le revêtement $\Vreg \surj W \qg \Vreg$ permet
  en effet d'obtenir la suite exacte
  \[ 1 \to P(W)=\pi_1(\Vreg) \to B(W)=\pi_1 (W \qg \Vreg)
  \stackrel{\pi}{\to} W \to 1 .\] Pour définir précisément $\pi$, on doit
  choisir une section $\widetilde{\CU}$ du point-base $\CU$ dans
  $\Vreg$. Il y a $|W|$ choix possibles, qui donnent des morphismes
  conjugués (cf. \cite[Rq.6.4]{BessisKPi1}).
\end{rqe}

Géométriquement, $s_i'$ est représenté par un chemin qui, vu depuis
$\CU$, passe sous $x_i$ mais pas sous $x_j$ pour $j\neq i$ :
{\shorthandoff{;:}%
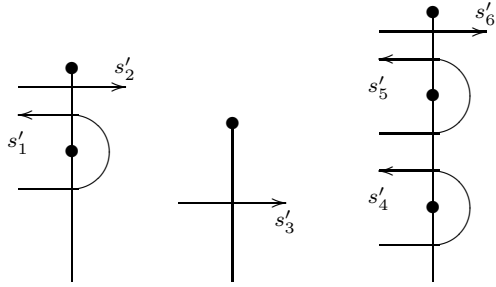
 \[\xy
<5pt,0pt>:<0pt,3.5pt>::
(0,0)="A", (0,-17)="A1", 
(-12,6)="B", (-12,-3)="B1", (-12,-17)="B2", 
(15,12)="C", (15,3)="C1", (15,-9)="C2", (15,-17)="C3", 
"B"*{\bullet},"B1"*{\bullet},"A"*{\bullet},"C"*{\bullet},"C1"*{\bullet},"C2"*{\bullet},
% (-15,6)*{_{x_2}},(-15,-3)*{_{x_1}},(-3,0)*{_{x_3}},(12,-9)*{_{x_4}},(12,3)*{_{x_5}},(12,12)*{_{x_6}},
(-16,-2)*{_{s_1'}},(-8,6)*{_{s_2'}},(4,-10.5)*{_{s_3'}},(11,-8)*{_{s_4'}},(11,4)*{_{s_5'}},(19,12)*{_{s_6'}},
(-11.5,-7);(-16,-7) **@{-},(-11.5,1);(-16,1) **@{-}, *\dir{>},
"B1"*\cir<14pt>{r^l},%flèche s_1'
(-16,4);(-8,4) **@{-}, *\dir{>},%flèche s_2'
(-4,-8.5);(4,-8.5) **@{-}, *\dir{>},%flèche s_3'
(15.5,-5);(11,-5) **@{-}, *\dir{>},(15.5,-13);(11,-13) **@{-},
"C2"*\cir<14pt>{r^l},%flèche s_4'
(15.5,7);(11,7) **@{-},*\dir{>},(15.5,-1);(11,-1) **@{-},
"C1"*\cir<14pt>{r^l},%flèche s_5'
(11,10);(19,10) **@{-}, *\dir{>},%flèche s_6'
"B";"B1" **@{-},
"B2";"B1" **@{-},
"A";"A1" **@{-},
"C";"C1" **@{-},
"C1";"C2" **@{-},
"C2";"C3" **@{-}
\endxy \]
}

\begin{rqe}
\label{rqzariski}
Si $y \in Y-\CK$, $\fact_B(y)$ est un $n$-uplet de «~réflexions tressées~»
(terminologie de Broué, cf. \cite{broumarou}), et détermine des générateurs
de la monodromie de $W \qg \Vreg=\BC^n - \CH$. Plus précisément, par
$\pi_1$-surjectivité de l'inclusion ${L_y - L_y\cap \CH } \ \inj \ {\BC^n -
  \CH}$ (cf. \cite[Thm. 2.5]{zariski}), les facteurs de $\fact_B(y)$
engendrent le groupe de tresses $B(W)$.
\end{rqe}

Dans le cas où le support de $\LL(y)$ est «~générique~» (parties réelles
distinctes), les $k$-uplets $\fact_B(y)$ et $\lbl(y)$ coïncident. Dans le
cas général, on peut toujours perturber $y$ en un $y'$ tel que
$\LL(y')=e^{-i\theta}\LL(y)$, avec $\theta>0$ assez petit pour que le
support de $\LL(y')$ soit générique et que
$\fact_B(y)=\fact_B(y')=\lbl(y')$. Les propriétés de $\lbl$ énoncées dans
\cite{BessisKPi1} s'adaptent ainsi aisément à l'application $\fact$ ; par
la suite (\ref{propfacto}, \ref{propcompat}, \ref{thmbij}), on rappelle ces
propriétés, en les reformulant dans ce nouveau cadre, et on ne démontre que
ce qui change de façon non triviale.

\medskip

Si $\LL(y)=\{0\}$ (avec multiplicité $n$), alors $y=0$ (d'après le théorème
\ref{thmLLBessis}). On note alors $\delta$ l'élément de $B(W)$ tel que
$\fact_B(0)=(\delta)$. Il est représenté par l'image dans $W \qg \Vreg$ du
chemin dans $\Vreg$ :
\begin{eqnarray*}
  [0,1] & \longrightarrow & V^{\reg}\\
  t & \longmapsto & v\exp(2i\pi t/h).
\end{eqnarray*}
où $v$ est tel que pour $i=1,\dots,n-1$, $f_i(v)=0$ (\ie $v\in
L_0$). Notons que $\delta$ est une racine $h$-ième du «~tour complet~» de
$P(W)$ («~full-twist~» noté $\bm{\pi}$ dans \cite[Not.2.3]{broumarou}).

On pose $c=\pi(\delta)$, image de $\delta$ dans $W$. Par construction, $c$
est $\eh$-régulier, donc est un élément de Coxeter. Les autres éléments de
Coxeter de $W$ (conjugués) sont obtenus pour d'autres choix du le morphisme
$\pi$ (cf. remarque \ref{rqpi}).

\medskip

Le lemme 6.14 dans \cite{BessisKPi1} permet de comparer les tunnels lorsqu'on
change de fibre $L_y$ :

\begin{lemme}[«~Règle d'Hurwitz~»]
  \label{reglehur}
  Soit $T=(y,z,L)$ un tunnel, qui représente un élément $s\in B(W)$. Soit
  $\Omega$ un voisinage connexe par arcs de $y$, tel que pour tout $y'\in
  \Omega$, $(y',z,L)$ soit encore un tunnel. Alors, pour tout $y'\in
  \Omega$, $(y',z,L)$ représente $s$.
\end{lemme}

On rassemble dans la proposition suivante quelques conséquences
utiles de ce lemme, adaptées de \cite[Lemmes 6.16, 6.17]{BessisKPi1} :

\begin{prop}
\label{propfacto}
  Soit $y\in Y$, $(x_1,\dots, x_k)$ le support ordonné de $\LL(y)$, et
  $(s_1,\dots,s_k)=\fact_B(y)$. Alors :
  \begin{enumerate}[(i)]
  \item $s_1\dots s_k=\delta$ et $\pi(s_1)\dots \pi(s_k)=c$ ;
  \item pour tout $i$, la longueur de $\pi(s_i)\in W$ est la multiplicité
    de $x_i$ dans $\LL(y)$.
  \end{enumerate}
  Par conséquent, $\fact(y)\in \Fact(c)$, et si $y \in \LL ^{-1}
  (E_{\lambda}^0)$, $\fact(y) \in \Fact_\lambda (c)$.
\end{prop}

Une autre conséquence est la compatibilité des actions d'Hurwitz et de
Galois :
\begin{prop}[d'après {\cite[Cor.6.18]{BessisKPi1}}]
  \label{propcompat}
  Soit $y\in Y-\CK$, et $\beta \in B_n$. Alors :
  \[ \fact(y \cdot \beta) = \fact(y)\cdot \beta \] où la première action
  est l'action de Galois définie en \ref{subpartLL}, et la seconde est
  l'action d'Hurwitz (cf. définition \ref{defhur}).
\end{prop}

\begin{proof}[Démonstration :]
  La preuve de \cite[Cor.6.18]{BessisKPi1} s'adapte au cas non générique
  (parties réelles non distinctes) sans difficultés. On verra plus loin
  (lemme \ref{lemcompat}) une généralisation de cette propriété de
  compatibilité à tous les éléments de $Y$.
\end{proof}

\section{\texorpdfstring{Etude de $\LL$ sur les strates $E_\lambda$}{Etude de LL sur les strates}}

\label{partLLstrat}

Ci-dessous on reformule \cite[Thm. 7.9]{BessisKPi1} en utilisant le produit
fibré $E_n \times_{\cp(n)} \Fact(c)$, où $\cp(n)$ désigne l'ensemble des
compositions de $n$. A toute configuration $X$ de $E_n$ on peut associer
une composition de $n$, constituée des multiplicités des points du support
de $X$ pris dans l'ordre $\lex$ sur $\BC$. D'autre part, toute
factorisation de $\Fact(c)$ détermine une composition de $n$, en
considérant les longueurs des facteurs dans l'ordre. Ces constructions
définissent deux applications : $\comp_1 : E_n \to \cp(n)$ et $\comp_2 :
\Fact(c) \to \cp(n)$. On pose :
\[ E_n \times_{\cp(n)} \Fact(c) := \{(X,\xi) \in E_n \times \Fact(c)
\tq \comp_1(X)=\comp_2(\xi) \} \]

\begin{thm}[d'après {\cite[Thm. 7.9]{BessisKPi1}}]
  \label{thmbij}
  L'application $\LL \times \fact$ :
  \[ \begin{array}{lcl}
    Y & \to & E_n  \times_{\cp(n)}  \Fact(c)\\
    y & \mapsto &(\LL(y) , \quad \fact(y))
  \end{array}
  \]
  est bijective.
\end{thm}

Donnons ici les grandes lignes de la démonstration (détails dans
\cite[Part.7]{BessisKPi1}). L'essentiel est de montrer la bijectivité de
$\LL\times \fact : Y-\CK \to \Red(c)$ (le reste se fait en dégénérant les
points réguliers). Cela revient à prouver que, pour $y \in Y-\CK$, le
morphisme de $B_n$-ensembles
\[ y\cdot B_n \xrightarrow{\fact} \fact(y)\cdot B_n \] est un isomorphisme
\cite[Thm. 7.4]{BessisKPi1}. Pour cela, on montre que ${|\fact(y)\cdot
  B_n|}={|y\cdot B_n|}=\frac{n!h^n}{|W|}$ en utilisant les deux propriétés
suivantes, prouvées au cas par cas \cite[Prop.7.5]{BessisKPi1} :
  \begin{enumerate}[(i)]
  \item l'action d'Hurwitz sur $\Red(c)$ est transitive ;
  \item $|\Red(c)|= \frac{n!h^n}{|W|}$.
  \end{enumerate}

\bigskip

Posons maintenant $Y_{\lambda} := \LL ^{-1} (E_{\lambda})$, et $Y_\lambda^0:=\LL
^{-1} (E_{\lambda}^0)$. Ainsi, $Y_{\lambda}$ est un fermé de Zariski
dans $Y\simeq \BC^{n-1}$, et est l'adhérence de $Y_{\lambda}^0$. Pour
tout $\lambda \vdash n$, $Y_\lambda = \bigsqcup_{\mu \leq \lambda}
Y_\mu ^0$. En particulier, $Y=\bigsqcup_{\mu \vdash n} Y_\mu ^0$.

\begin{thm}
  \label{thmrevetement}
  Pour tout $\lambda \vdash n$, la restriction de $\LL$ :
  \[ \LL_{\lambda} :Y_{\lambda}^0 \ \surj \ E_{\lambda}^0 \] est un
  revêtement non ramifié.
\end{thm}

\begin{rqe}
  La définition de «~revêtement~» est bien sûr à prendre ici au sens large
  : comme $Y_{\lambda}^0$ n'est pas nécessairement connexe par arcs, le
  théorème signifie que l'application $\LL_\lambda$ est un revêtement
  connexe par arcs sur chacune des composantes connexes par arcs de
  $Y_{\lambda}^0$.
\end{rqe}

\begin{proof}[Démonstration :]
  Soit $X_0 \in E_{\lambda}^0$. Supposons dans un premier temps que les
  parties réelles des éléments du support de $X_0$ sont distinctes (on dira
  que $X_0$ est relativement générique). Soit $\Omega$ un ouvert connexe
  par arcs de $E_{\lambda}^0$, contenant $X_0$, et de diamètre assez petit
  pour que tous les éléments de $X_0$ soient relativement génériques.

  Soit $\mu=\comp_1(X_0)$ la composition de $n$ associée à $X_0$. Il est
  clair que pour tout $X \in \Omega$, la composition associée à $X$ est
  encore $\mu$. Posons $F_\mu:=\Fact_{\mu} (c)=\{\xi \in \Fact(c) \tq
  \comp_2(\xi)=\mu \}$. Le théorème \ref{thmbij} implique que l'application
  \[ \LL_{\lambda}^{-1}(\Omega) \xrightarrow{\LL_{\lambda}\times \fact}
  \Omega \times F_\mu \] est une bijection.

    Si $X_0$ n'est pas relativement générique, on peut refaire toutes les
  constructions précédentes en tournant légèrement la direction de la
  verticale dans le sens trigonométrique direct. On doit modifier alors la
  définition de l'application $\fact$ (sauf pour les $y$ de la fibre de
  $X_0$), mais les factorisations construites restent toujours dans
  $\Fact_{\mu} (c)$ pour $\mu = \comp_1 (X_0)$.

  Pour conclure, il suffit de remarquer que toutes les fibres $F_\mu$
  construites sont en bijection. En effet, les compositions $\mu$ associées
  aux éléments de $Y_\lambda ^0$ correspondent nécessairement à la
  partition $\lambda$, et si deux compositions $\mu$ et $\mu'$ sont égales
  à permutation des parts près, alors $\Fact_\mu (c)$ et $\Fact_{\mu'}(c)$
  sont en bijection (utiliser l'action d'Hurwitz par une tresse adaptée).
\end{proof}

Soit $\lambda$ une partition de $n$, et $X,X'\in E_\lambda ^0$. On a vu en
partie \ref{subpartEn} que toute classe d'homotopie de chemin de $X$ vers
$X'$ dans $E_\lambda ^0$ définit un élément de $B_p$ où $p=\# \lambda$. On
en déduit en particulier une action de Galois de $B_p$ sur chacune des
fibres de $\LL$ au-dessus de $E_{\lambda}^0$, et on a une propriété de
compatibilité des actions plus générale :

\begin{lemme}[Compatibilité des actions de Galois et d'Hurwitz sur
  les strates]
  \label{lemcompat}
  Soit $\lambda \vdash n$, avec $p=\#\lambda$.
  \begin{enumerate}[(i)]
  \item Soient $y,y' \in Y_\lambda^0$, reliés par un chemin $\gamma$ dans
    $Y_\lambda ^0$. Notons $X=\LL(y)$, $X'=\LL(y')$, et $\beta$ la tresse
    de $B_p$ représentée par l'image de $\gamma$ par $\LL$. Alors :
    $\fact(y')=\fact(y) \cdot \beta$.
  \item Soient $\beta \in B_p$, et $y,y' \in Y_\lambda ^0$ tels que
    $\fact(y')=\fact(y) \cdot \beta$. Notons $X= \LL(y)$, $X'=\LL(y')$ et $
    \widetilde{\beta}$ l'unique relevé (par $\LL$) d'origine $y$ de la
    tresse $\beta$ vue comme classe de chemin (dans $E_\lambda^0$) de $X$
    vers $X'$. Alors $y'$ est le point d'arrivée de $\widetilde{\beta}$.
  \end{enumerate}
\end{lemme}

\begin{proof}[Démonstration :]
  \begin{enumerate}[(i)]
  \item C'est essentiellement la même preuve que pour
    \cite[Cor.6.18]{BessisKPi1}. Il suffit de considérer le cas
    $\beta=\bm{\sigma_i}$ ($i$-ème tresse élémentaire de $B_p$), pour $i
    \in \{1,\dots, p-1\}$. En vertu du théorème \ref{thmrevetement}, on
    peut déplacer $X,X'$ (sans changer l'ordre des points des
    configurations), et $y,y'$ (sans modifier $\fact(y)$ et $\fact(y')$),
    de sorte que la tresse $\bm{\sigma_i}$ soit représentée par le chemin
    suivant (en pointillés) : {\shorthandoff{;:}%
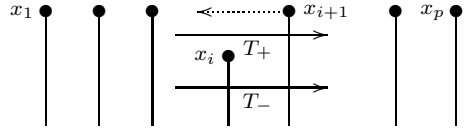
      \[\xy
      (-10,0)="1", (-3,0)="2", (4,0)="3", 
      (14,-6)="4", (22,0)="5", (36,0)="6", (44,0)="7",
      (-10,-15)="11", (-3,-15)="22", (4,-15)="33", 
      (14,-15)="44", (22,-15)="55", (36,-15)="66", (44,-15)="77",
      "1";"11" **@{-}, "2";"22" **@{-}, "3";"33" **@{-}, "4";"44" **@{-},
      "5";"55" **@{-}, "6";"66" **@{-}, "7";"77" **@{-},
      "1"*{\bullet},"2"*{\bullet},"3"*{\bullet},"4"*{\bullet},
      "5"*{\bullet},"6"*{\bullet},"7"*{\bullet},
      (11,-6)*{_{x_i}},(27,0)*{_{x_{i+1}}},
      (-13,0)*{_{x_1}},(41,0)*{_{x_{p}}},
      "5";(10,0) **@{.}, *\dir{>},
      (7,-3);(27,-3) **@{-}, *\dir{>},
      (7,-10);(27,-10) **@{-}, *\dir{>},
      (18,-12)*{_{T_-}},
      (18,-5)*{_{T_+}}
      \endxy \] } Considérons les tunnels $T_+$ et $T_-$ représentés sur le
    schéma. Notons $(w_1,\dots, w_p)=\fact(y)$ et $(w_1',\dots,
    w_p')=\fact(y')$. Alors $T_+$ représente $w_{i+1}$ dans $L_y$, et
    $w_i'$ dans $L_{y'}$. D'où, d'après le lemme \ref{reglehur},
    $w_i'=w_{i+1}$. De même, avec $T_-$, on obtient $w_i'w_{i+1}'=w_i
    w_{i+1}$, et on a vérifié que $\fact(y')=\fact(y) \cdot \bm{\sigma_i}$.
  \item On reprend les notations de l'énoncé. Notons $y''$ le point
    d'arrivée de $\beta$. Alors, d'une part on a
    $\LL(y'')=X'=\LL(y)$. D'autre part, en appliquant le point (i) à $y$,
    $y'$, et $\widetilde{\beta}$, on obtient : $\fact(y'')=\fact(y)\cdot
    \beta=\fact (y')$. D'où, par le théorème \ref{thmbij}, $y''=y'$.
  \end{enumerate}
\end{proof}

Le théorème suivant est une conséquence directe du lemme :

\begin{thm}
  \label{thmcomposantes}
  Soit $\lambda \vdash n$, et $p=\#\lambda$. L'application $ Y_\lambda^0
  \xrightarrow{\fact} \Fact_\lambda(c)$ induit une bijection entre
  l'ensemble des composantes connexes par arcs de $Y_\lambda^0$ et
  l'ensemble des orbites d'Hurwitz de $\Fact_\lambda(c)$ sous $B_p$.

  Autrement dit, si $ Y_{\lambda}^0 = \bigsqcup_{i} Y_{\lambda,i}^0$ est la
  décomposition de $Y_{\lambda}^0$ en ses composantes connexes par arcs,
  alors
  \[ \Fact_{\lambda}(c)=\fact(Y_{\lambda}^0) = \bigsqcup_{i}
  \fact(Y_{\lambda,i}^0) \] est la décomposition de $\Fact_{\lambda}(c)$ en
  orbites d'Hurwitz sous $B_p$.
\end{thm}

\section{Stratification de $\CH$ et éléments de Coxeter paraboliques}
\label{partecp}

Dans cette partie on étudie la géométrie de $\CH$, afin de pouvoir
déterminer en partie \ref{partirred} les composantes connexes par arcs de
$Y_\lambda ^0$ lorsque $\lambda$ est une partition primitive, en appliquant
le théorème \ref{thmcomposantes} ci-dessus.

On rappelle qu'on a fixé un élément de Coxeter $c=\pi(\delta)$.

\subsection{Stratification de $V$}
~ 

L'arrangement d'hyperplans associé à $W$ est noté $\CA$. On
considère la stratification de $W$ par les \emph{plats}, \ie
le treillis d'intersection :

\[ \CL:=\CL(\CA)=\left\{\bigcap_{H\in \CB} H \tq \CB \subseteq \CA \right\} .\]

Pour $L \in \CL$, on pose 
\[L^0:=L - \bigcup_{L'\in \CL, L' \subsetneq L} L'.\]
On obtient la stratification ouverte de $V$ associée à
$\CL$. Pour tout $L\in \CL$, $L$ est l'adhérence de $L^0$.

La stratification de $V$ par les plats correspond à la stratification du
groupe $W$ par ses sous-groupes paraboliques : on rappelle le théorème de
Steinberg \cite[Thm. 1.5]{steinberg} et ses conséquences.

\begin{thm}[Steinberg]
\label{thmsteinberg}
  Si $L$ est un plat, le groupe 
  \[W_L:=\{w\in W \tq \forall x \in L, wx=x\} \]
  est encore un groupe de réflexions, appelé sous-groupe parabolique
  associé à $L$. De plus :
  \begin{enumerate}[(i)]
  \item L'application $L \mapsto W_L$ est une bijection de $\CL$ vers
    l'ensemble des sous-groupes paraboliques ; sa réciproque est
    \[ G \mapsto V^{G}=\{x\in V \tq \forall w \in G, wx=x \}=\bigcap_{r
      \in \CR\cap G} H_r \ .\]
  \item Le rang de $W_L$ est égal à la codimension de $L$.
  \item Soit $v \in V$. Notons $V_v:=\displaystyle{\bigcap_{H\in \CA, v\in H} H}$
    et $W_v:= \{w\in W \tq wv=v\}$.\\ Alors, pour $L\in \CL$ :
    \[ v \in L^0 \ \ssi \  V_v=L \ \ssi \ W_v=W_L .\]
  \end{enumerate}
\end{thm}

\subsection{Éléments de Coxeter paraboliques}
~
\label{subpartcoxparab}

Notons $f :\CH \to \NCP_W(c)$ l'application qui à $(y,x) \in \CH$ associe
le facteur de $\fact(y)$ correspondant à $x$ (\ie défini par un tunnel
élémentaire passant sous le point $x$ de $\LL(y)$). Ainsi,
$(f(y,x_1),\dots, f(y,x_p))=\fact(y)$ lorsque le support ordonné de
$\LL(y)$ est $(x_1,\dots, x_p)$ .
 
On reformule ci-dessous un lemme fondamental de \cite{BessisKPi1} :

\begin{lemme}[d'après {\cite[Lemme 7.3]{BessisKPi1}}]
  \label{lemcoxpar}
  Soit $y\in Y$. Soient $x \in \LL(y)$, de multiplicité $p$, et $w=f(y,x)$.

  Alors il existe une préimage $v\in V$
  de $(y,x)\in W \qg V$ telle que $w$ soit un élément
  de Coxeter dans le sous-groupe parabolique $W_{v}$. De plus :
  \[ p=\ell(w)=\rg W_{v}=\dim V/V_{v} = \codim \Ker (w-1). \]
\end{lemme}

Par conséquent, en utilisant la surjectivité de l'application $\fact$, on
peut déduire de ce lemme que tout diviseur $w$ de $c$ est un élément de
Coxeter d'un sous-groupe parabolique, groupe que l'on peut déterminer à
l'aide d'une factorisation qui contient $w$. D'où la proposition-définition
suivante :

\begin{prop}
  \label{propcoxpar}
  Soit $W$ un groupe de réflexions complexe bien engendré, et $w \in
  W$. Les propriétés suivantes sont équivalentes :
  \begin{enumerate}[(i)]
  \item $w$ est un élément de Coxeter d'un sous-groupe parabolique de $W$ ;
  \item il existe un élément de Coxeter $c_w $ de $W$, tel que $w \< c_w$ ;
  \item $w$ est conjugué à un élément de $\NCP_W(c)$.
  \end{enumerate}
  On dit alors que $w$ est un \emph{élément de Coxeter parabolique}.
\end{prop}

\begin{rqe}
Dans le cas d'un groupe de Coxeter fini, cette propriété est bien connue et
démontrée de manière uniforme, cf. \cite[Lemme 1.4.3]{Bessisdual}.
\end{rqe}

\begin{proof}[Démonstration :]
  (ii) $\ssi$ (iii) provient directement de la théorie de Springer (les
  éléments de Coxeter de $W$ forment une seule classe de conjugaison) et de la
  propriété $\NCP_W(a c a ^{-1}) = a \NCP_W(c) a ^{-1}$.

  (i) $\Rightarrow$ (ii) : soit $G$ un sous-groupe parabolique (non
  trivial), et $w$ un élément de Coxeter de $G$. Soit $L=V^G$ et $v \in
  L^0$, de sorte que $W_v=G$. Notons $(y,x)=\bar{v}$ et $w_0=f(y,x) \in
  \NCP_W(c)$. D'après le lemme \ref{lemcoxpar}, il existe $v_0\in V$ tel
  que $(y,x)=\overline{v_0}$ et que $w_0$ soit un élément de Coxeter de
  $W_{v_0}$. Comme $v_0$ et $v$ sont dans la même orbite sous $W$, $G=W_v$
  est conjugué à $W_{v_0}$, donc tous leurs éléments de Coxeter sont
  conjugués dans $W$. En particulier $w$ est conjugué à $w_0$.

  (iii) $\Rightarrow$ (i) : il suffit de montrer l'implication pour $w \in
  \NCP_W(c)$, puisque la propriété (i) est invariante par conjugaison. Si
  $w \<c$, la surjectivité de $\fact$ (théorème \ref{thmbij}) donne
  l'existence de $(y,x)\in W \qg V $ tel que $f(y,x)=w$. Le lemme
  \ref{lemcoxpar} permet alors de conclure.
\end{proof}

Comme dans le cas d'un groupe de Coxeter \cite[Cor.1.6.2]{Bessisdual}, on
peut retrouver, à partir d'un élément de Coxeter parabolique, le
sous-groupe parabolique associé :

\begin{prop}
  \label{propcoxpar2}
  Soit $w$ un élément de Coxeter parabolique, et $W_w$ le sous-groupe
  parabolique fixateur du plat $\Ker(w-1)$. Alors :
  \begin{enumerate}[(i)]
  \item le groupe $W_w$ est l'unique sous-groupe parabolique duquel $w$ est
    un élément de Coxeter ;
  \item si $(r_1,\dots,r_k)\in \Red (w)$, alors $\left< r_1,\dots, r_k
    \right>=W_w$.
  \end{enumerate}
\end{prop}

\begin{proof}[Démonstration :]
  (i) Soit $G$ un sous-groupe parabolique tel que $w$ soit un élément de
  Coxeter de $G$. Soit $L$ le plat associé à $G$, \ie $L=V^{G}$. Alors
  comme $w\in G$, on a : $L \subseteq \Ker(w-1)$. Or $\codim L =\rg G$ par
  théorème \ref{thmsteinberg}, et $\rg G=\ell_{G} (w)$ (où $\ell_G$ désigne la
  longueur relativement aux réflexions de $G$) car $w$ est un élément de
  Coxeter. D'autre part $\ell(w)=\codim \Ker(w-1)$ par la proposition
  \ref{proplongueur}. Donc pour conclure il suffit d'utiliser le résultat
  suivant :
  \[ \mathrm{Si}\ G\ \text{est\ un\ sous-groupe\ parabolique,\ alors\ :\ }
  \forall g \in \NCP_W,\ g\in G \Rightarrow \ell_G(g)=\ell(g). \] Pour cela, on
  va vérifier que pour $r\in \CR$, si $r\< g$, alors $r \in G$. Si $G=W_L$,
  cela revient à montrer que $L \subseteq \Ker(r-1)$. Or, comme $g\in G$,
  on a $L\subseteq \Ker(g-1)$, et d'après la proposition \ref{propordre},
  $\Ker(g-1) \subseteq \Ker(r-1)$.

  (ii) Il suffit de le démontrer dans le cas où $w$ est un élément de
  Coxeter $c$ d'un groupe de réflexions (bien engendré) irréductible
  $W$. D'après le théorème \ref{thmbij}, toute décomposition réduite
  $(r_1,\dots, r_n)$ de $c$ provient d'une factorisation de $\delta$ en
  générateurs de la monodromie, \ie de $\fact_B (y)=(s_1,\dots, s_n)$, avec
  $y \in Y-\CK$ et $\pi(s_i)=r_i$. Or on sait qu'alors $s_1,\dots, s_n$
  engendrent $B(W)$ (cf. remarque \ref{rqzariski}), donc $r_1,\dots, r_n$
  engendrent $W$.
\end{proof}

Dans le cas présent, on a fixé un élément de Coxeter $c$ de $W$, et on ne
va considérer que les éléments paraboliques qui sont dans $\NCP_W$.  On
n'obtient donc pas tous les sous-groupes paraboliques, mais seulement les
\emph{sous-groupes paraboliques «~non croisés~»}, \ie ceux qui possèdent un
élément de Coxeter qui divise $c$. Cependant d'après la proposition
\ref{propcoxpar2}, quitte à les conjuguer, on obtient tous les sous-groupes
paraboliques :

\begin{prop}
\label{propsgpnc}
Soit $W$ un groupe de réflexions complexe bien engendré, et $c$ un élément
de Coxeter fixé de $W$. Soit $W_0$ un sous-groupe parabolique de $W$.

Alors $W_0$ est conjugué à un sous-groupe parabolique «~non croisé~» de
$W$, \ie un sous-groupe de la forme $W_L$ où $L=\Ker(w-1)$ et $w \< c$.
\end{prop}

Ce un résultat non trivial utilise ainsi de façon essentielle le lemme
\ref{lemcoxpar}. On peut dire que les sous-groupes paraboliques non croisés
jouent ici le rôle des sous-diagrammes (ou des paraboliques standards) de
la théorie de Coxeter (voir aussi \cite[p.3]{broumarou1} sur les diagrammes
et sous-diagrammes pour les groupes complexes).

\begin{rqe}
  Dans un sous-groupe parabolique non croisé, il y a unicité de l'élément
  de Coxeter divisant $c$. En effet, il faut trouver un diviseur de $c$
  dont le plat associé est donné, et la solution est unique par le théorème
  de Brady-Watt (\ref{thmBW}). L'ensemble des sous-groupes paraboliques non
  croisés, ordonné par inclusion, est donc isomorphe au treillis
  $(\NCP_W(c), \<)$.
\end{rqe}

\subsection{\texorpdfstring{Stratification de $W \qg V$}{Stratification de
    W \textbackslash V}}
\label{subpartstratquotient}
~

Le groupe $W$ agit sur $\CA$, donc sur $\CL$. On peut ainsi définir des
orbites de plats, qui forment une stratification notée $\Lb$ de $W \qg V$. 

Soit $p$ la projection $V \rightarrow W \qg V$,
$v\mapsto \bar{v}=W\cdot v$. On a : $\Lb= W \qg \CL =(p(L))_{L \in \CL}=
(W\cdot L)_{L \in \CL} $.

\medskip

On notera par la suite les strates de $\Lb$ par la lettre $\Lambda$. Pour $
\Lambda \in \Lb $ , posons :
\[\Lambda^0:=\Lambda - \bigcup_{\Lambda' \in \Lb, \Lambda'
  \subsetneq \Lambda} \Lambda' \]

Si $\Lambda=W \cdot L$, alors $\Lambda^0=W \cdot L^0$. Les ouverts
$\Lambda^0$, pour $\Lambda\in \Lb$, forment la stratification ouverte de $W
\qg V$ associée à $\Lb$, appelée \emph{stratification
  discriminante}. Notons que $(W \qg V)^0=W \qg V - \CH = W \qg \Vreg$.

\bigskip

Les strates de $\Lb$ correspondent aux classes de conjugaison de
sous-groupes paraboliques, puisque $W_{w\cdot L} = w W_L w ^{-1}$. On peut
également les associer aux classes de conjugaison d'éléments de Coxeter
paraboliques. En effet, considérons l'application $F: W\to \CL$, $w \mapsto
\Ker(w-1)$. Si $w$ et $w'$ sont conjugués, alors $F(w)$ et $F(w')$ sont
dans la même orbite sous $W$, donc $F$ induit une application $\bar{F}$ de
l'ensemble des classes de conjugaison de $W$ vers $\Lb$.

\begin{prop}
  L'application $\bar{F}$ définie ci-dessus induit une bijection entre
  l'ensemble $\Lb$ des strates de $W \qg V$ et :
  \begin{itemize}
  \item l'ensemble des classes de conjugaison d'éléments de Coxeter
    paraboliques ;
  \item l'ensemble des classes de conjugaison d'éléments de $\NCP_W$.
  \end{itemize}
\end{prop}

\begin{proof}[Démonstration :]
  En utilisant les propositions \ref{propcoxpar} et \ref{propcoxpar2}, le
  premier point est clair : deux éléments de Coxeter d'un même sous-groupe
  parabolique sont conjugués d'après la théorie de Springer, donc deux éléments
  de Coxeter paraboliques $w_1$ et $w_2$, associés à deux sous-groupes
  paraboliques $W_1$ et $W_2$, sont conjugués si et seulement si les
  groupes $W_1$ et $W_2$ sont conjugués.

  Le second point est direct en utilisant la proposition \ref{propcoxpar}.
\end{proof}

\begin{defn}
  \label{deftype}
  Soit $\Lambda$ une strate de $\Lb$, $w \in \NCP_W$, et  $\bar{v} \in W \qg
  V$. On dit que :
  \begin{itemize}
  \item «~$w$ est de \emph{type} $\Lambda$~» si la classe de conjugaison de $w$
    correspond à $\Lambda$ par la bijection ci-dessus, \ie si $\Lambda = W \cdot
    \Ker(w-1)$.
  \item «~la \emph{strate} de $\bar{v}$ est $\Lambda$~» si $\Lambda$ est la
    strate \emph{minimale} de $\Lb$ contenant $\bar{v}$, \ie si $\bar{v}
    \in \Lambda^0$, ou encore $\Lambda=W \cdot V_v$.
  \end{itemize}
\end{defn}

On peut ainsi reformuler le lemme \ref{lemcoxpar} :

\begin{lemme}
  \label{lemtype}
  Soit $(y,x) \in \CH$. Alors le type de l'élément de Coxeter parabolique
  $f(y,x)$ est la strate du point $(y,x)$.
\end{lemme}

Une conséquence du lemme est que si $(y,x)$ est dans une
strate ouverte de dimension $n-k$, alors $\ell(f(y,x))=k$ et la
multiplicité de $x$ dans $\LL(y)$ est $k$. Soit $\Lb_k$ l'ensemble des
strates fermées de dimension exactement $n-k$, et posons :

\[ \CH_k:= \bigcup_{\Lambda \in \Lb_k} \Lambda  \]

Ainsi $W \qg V= \CH_0 \supsetneq \CH_1=\CH \supsetneq \CH_2 \supsetneq
\dots \supsetneq \CH_n = \{0 \}$. 

\bigskip

Désormais on suppose $k\geq 1$. Soit la projection $\varphi : \CH \to Y,\
(y,x)\mapsto y$. Notons $Y_k := \varphi(\CH_k)$, et $\alpha_k := k^1
1^{n-k} \vdash n$. D'après le lemme \ref{lemcoxpar}, pour $(y,x)\in \CH$,
la longueur du facteur $f(y,x)$ est donnée par la codimension de la strate
de $(y,x)$, d'où la propriété suivante :

\begin{lemme}
  Soit $(y,x)\in W \qg V$. Alors  $(y,x)$ est dans $\CH_k$ si et
  seulement si la multiplicité de $x$ dans $\LL(y)$ est supérieur ou égale
  à $k$.

  Par conséquent : $Y_k= \LL ^{-1} (E_{\alpha_k})$.

\end{lemme}

La partie $Y_k$ est donc ce que l'on avait noté $Y_\lambda$ dans la partie 5,
dans le cas où $\lambda$ est la partition primitive $\alpha_k$. En
particulier, comme $\CK=\LL ^{-1} (E_{\alpha_2})$, on obtient $Y_2=\CK$.

\section{Composantes connexes par arcs de $Y_k^0$}
\label{partirred}

D'après le théorème \ref{thmcomposantes}, pour étudier les orbites
d'Hurwitz primitives, il est intéressant d'identifier les composantes
connexes par arcs de $Y_k^0$. Pour cela, on commence par déterminer les
composantes irréductibles de $Y_k$.

\subsection{Composantes irréductibles de $\CH_k$ et de $Y_k$}
~

On va d'abord montrer que les strates de $\Lb_k$ sont les composantes
irréductibles de $\CH_k$.

\begin{lemme}
  \label{lemfini}
  Les morphismes $p: V \surj W \qg V$ et $\varphi:\CH \rightarrow Y$ sont
  finis. Par conséquent ils sont fermés, pour la topologie de Zariski.
\end{lemme}

\begin{proof}[Démonstration :] 
  Pour $p$, c'est le théorème de Chevalley : $\CO_V \simeq \BC[v_1,\dots,
  v_n]$, $\CO_{W \qg V} \simeq \BC[f_1,\dots, f_n]$, et $\CO_V$ est un
  $\CO_{W \qg V}$-module libre de rang $|W|$.

  Pour $\varphi$, on a $\CO_Y=\BC[f_1,\dots, f_n]$ et $\CO_\CH =
  \BC[f_1,\dots, f_n] / (\D)$, avec $\D=f_n^n+a_2 f_n^{n-2} + \dots + a_n$
  et $a_i \in \BC[f_1,\dots, f_{n-1}]$. En particulier, $f_n$ est entier sur
  $\varphi^*(\CO_Y)$. Donc $\varphi^*$ fait de $\CO_\CH$ un $\CO_Y$-module
  libre de rang $n$.
\end{proof}

Par conséquent, les strates de $\Lb$, images par $p$ des plats dans $V$,
sont fermées, et irréductibles (car l'image d'un irréductible par un
morphisme algébrique est irréductible). D'où :

\begin{cor}
  Pour tout $k \geq 1$, les composantes irréductibles du fermé $\CH_k$ sont
  les strates $\Lambda$ de $\Lb_k$.
\end{cor}

Par $\varphi$, continue et fermée, on peut envoyer ces strates dans
$Y_k=\varphi(\CH_k)$. On a ainsi $Y_k=\bigcup_{\Lambda \in \Lb_k}
\varphi(\Lambda)$, avec les $\varphi(\Lambda)$ fermés irréductibles dans
$Y$.

\medskip

Commençons par le cas $k=1$, qui correspond aux classes de conjugaison
de réflexions de $\NCP_W$. En utilisant le théorème \ref{thmbij}, on
obtient le résultat suivant.

\begin{prop}
  Soit $\Lambda \in \Lb_1$. Alors, pour tout $y \in Y - \CK$, au moins un
  des facteurs de $\fact(y)$ a pour classe de conjugaison $\Lambda$.

  Par conséquent, pour tout $\Lambda \in \Lb_1$, on a : $\varphi(\Lambda) = Y$.

\end{prop}

\begin{proof}[Démonstration :]

  D'après la surjectivité de l'application $\fact$, toute réflexion de
  $\NCP_W$ apparaît dans une factorisation $\fact(y)$. Donc, par
  transitivité de l'action d'Hurwitz sur $\Red(c)$, si $\xi \in
  \fact(Y-\CK)$, alors toutes les classes de conjugaison de réflexions de
  $\NCP_W$ apparaissent dans $\xi$. En effet, l'ensemble des classes de
  conjugaison de réflexions faisant partie d'une décomposition réduite est
  invariant par l'action d'Hurwitz.

  Soit $y \in Y$. Quitte à désingulariser, on peut trouver $y' \in Y-\CK$
  tel que $\fact(y')$ soit un raffinement de $\fact(y)$. Soit $\Lambda \in
  \Lb_1$, alors $\fact(y')$ contient un facteur de type $\Lambda$, donc il
  existe $x$ dans $\LL(y)$ tel que $f(y,x)$ soit multiple (pour $\<$) d'une
  réflexion de type $\Lambda$, d'où $(y,x) \in \Lambda$, et $y \in
  \varphi(\Lambda)$.
\end{proof}

Désormais on suppose $k \geq 2$.

\begin{prop}
  Les $\varphi(\Lambda)$, pour $\Lambda \in \Lb_k$, sont distincts deux à deux,
  et sont les composantes irréductibles de $Y_k$.
\end{prop}

\begin{proof}[Démonstration :]
  D'après la discussion précédente, il suffit de montrer que si $\Lambda,
  \Lambda' \in \Lb_k$, avec $\Lambda \neq \Lambda'$, alors $\varphi(\Lambda)
  \nsubseteq \varphi(\Lambda')$.
  
  La strate ouverte $\Lambda ^0$ correspond à une classe de conjugaison de
  sous-groupe parabolique de rang $k$. Soit $c_\Lambda$ un élément de Coxeter
  parabolique, divisant $c$, de type $\Lambda$. Complétons avec $n-k$ réflexions
  pour obtenir une factorisation complète : $\xi=(c_\Lambda,
  s_{k+1},\dots,s_n)$. Alors, par \ref{thmbij}, il existe $y$ dans $Y$,
  tel que $\fact(y)=\xi$, et que $\LL(y)$ ait $n-k+1$ points
  distincts. Soit $x$ le point multiple dans $\LL(y)$ ; alors $(y,x)\in
  \Lambda^0$ puisque l'élément de Coxeter parabolique associé à $x$ est de
  type $\Lambda$. D'où : $y\in \varphi(\Lambda)$.

  Supposons que $y\in \varphi(\Lambda')$ ; alors il existe $x'$ tel que
  $(y,x')\in \Lambda'$. Donc, dans $\fact(y)=\xi$, on doit trouver un
  élément de type $\Lambda''\subseteq \Lambda'$, de longueur supérieure ou
  égale à $k$ ; or, dans $\fact(y)$, seul $c_\Lambda$ convient, et il est
  de type $\Lambda \nsubseteq \Lambda'$. D'où $y \notin \varphi(\Lambda')$,
  ce qui conclut la preuve.
\end{proof}

\subsection{Connexité par arcs}
\label{subpartcpa}
~

Notons comme plus haut $Y_k^0:= \LL ^{-1} (E_{\alpha_k}^0)$, et pour $\Lambda \in
\Lb_k$, $\varphi(\Lambda)^0:=\varphi(\Lambda)\cap Y_k^0$. 

Pour $\Lambda \in \Lb_k$, notons $\Fact_{\alpha_k}^{\Lambda}(c)$ les factorisations
primitives de $c$, de forme $\alpha_k$, et dont l'élément long est de type
$\Lambda$. Alors, en vertu du lemme \ref{lemtype}, on a :
\[ \varphi(\Lambda)^0 = \fact ^{-1} (\Fact_{\alpha_k}^{\Lambda}(c)). \]

Ainsi $Y_k ^0=\bigsqcup_{\Lambda \in \Lb_k} \varphi(\Lambda)^0$. D'après la
surjectivité de $\fact$, $\varphi(\Lambda)^0$ est un ouvert (de Zariski) non
vide de $\varphi(\Lambda)$. 

\begin{prop}
\label{propcpa}
Pour tout $\Lambda \in \Lb_k$, $\varphi(\Lambda)^0$ est connexe par arcs.
\end{prop}

\begin{proof}[Démonstration :]
  C'est un fait général, mais non trivial, que tout ouvert de Zariski d'une
  variété algébrique complexe irréductible est connexe par arcs. Dans le
  cas présent on peut cependant donner une démonstration explicite. Soit
  $\Lambda \in \Lb_k$, et soit $L \in \CL$ tel que $p(L)=\Lambda$. Notons $\Omega=
  L\cap (\varphi \circ p)^{-1}(\varphi(\Lambda)^0)$. Comme $p: V \to W \qg V$ et
  $\varphi : W \qg V \to Y$ sont des morphismes algébriques, $\Omega$ est
  un ouvert de Zariski de $L$. Comme $L$ est un espace vectoriel sur $\BC$,
  la connexité par arcs de $\Omega$ est alors évidente. D'autre part,
  $\varphi(\Lambda)^0=\varphi \circ p (\Omega)$, avec $\varphi \circ p$
  continue, donc $\varphi(\Lambda)^0$ est connexe par arcs.
\end{proof}

En utilisant le théorème \ref{thmcomposantes}, on va en déduire aisément
que les $\varphi(\Lambda)^0$, pour $\Lambda \in \Lb_k$, sont bien les
composantes connexes par arcs de $Y_k^0$.

\section{ Forte conjugaison, et cas des réflexions}
\label{partrefl}

\begin{thm}
  \label{thmfconj2}
  Soit $k \in \{2,\dots , n \}$. Alors :
  \begin{itemize}
  \item les parties $\varphi(\Lambda)^0$, pour $\Lambda \in
    \Lb_k$, sont les composantes connexes par arcs de $Y_k^0$ ;
  \item les ensembles
    $\fact(\varphi(\Lambda)^0)=\Fact_{\alpha_k}^{\Lambda}(c)$ sont les
    orbites d'Hurwitz sous $B_{n-k+1}$ de
    $\fact(Y_k^0)=\Fact_{\alpha_k}(c)$.
  \end{itemize}
  Par conséquent, deux éléments de $\NCP_W$ de longueur $k$ sont fortement
  conjugués si et seulement s'ils sont conjugués.
\end{thm}

\begin{proof}[Démonstration :]
  Soient $y\in \varphi(\Lambda)^0$ et $y'\in \varphi(\Lambda')^0$, avec
  $\Lambda,\Lambda' \in \Lb_k$. Si $y$ et $y'$ sont reliés par un chemin
  dans $Y_k^0$, alors, par le théorème \ref{thmcomposantes}, $\fact(y)$ et
  $\fact(y')$ sont dans la même orbite d'Hurwitz, donc leurs éléments longs
  sont conjugués, \ie $\Lambda=\Lambda'$ (par lemme \ref{lemtype}). Donc la
  proposition \ref{propcpa} implique que les $\varphi(\Lambda)^0$ sont les
  composantes connexes par arcs de $Y_k^0$.  Les orbites d'Hurwitz de
  $\fact(Y_k^0)$ sont alors directement données par le théorème
  \ref{thmcomposantes}.

  Enfin, la propriété de forte conjugaison vient du fait qu'une
  factorisation de forme $\alpha_k$ est dans $\fact(\varphi(\Lambda)^0)$ si et
  seulement si son facteur long est de type $\Lambda$.
\end{proof}

Pour conclure la preuve du théorème \ref{thmfconj}, il reste à déterminer
les classes de conjugaison forte de réflexions :

\begin{thm}
  \label{thmfconjrefl}
  Soient $r,r'$ deux réflexions de $\NCP_W$. Si $r$ et $r'$ sont
  conjuguées, alors $r$ et $r'$ sont fortement conjuguées dans $\NCP_W$.
\end{thm}

\begin{rqe}
  Cette propriété apporte une précision intéressante concernant l'action
  d'Hurwitz de $B_n$ sur $\Red(c)$ : si $r_1, r_1' \in \NCP_W$ sont deux
  réflexions conjuguées, et si $(r_1,r_2,\dots, r_n)$ et $(r_1',r_2',\dots,
  r_n')$ sont deux décompositions réduites de $c$, alors il existe une
  tresse de $B_n$, \emph{pure par rapport au premier brin}, qui transforme
  l'une en l'autre (cf. remarque \ref{rqrefl} sur le lien entre forte
  conjugaison et action d'Hurwitz).
\end{rqe}

\begin{proof}[Démonstration :]
  Considérons $\varphi : \CH \to Y$, $(y,x) \mapsto y$. Notons $\CH' :=
  \varphi ^{-1} (Y- \CK)$. Alors la restriction $\varphi ' : \CH' \to Y -
  \CK$ est un revêtement non ramifié à $n$ feuillets (continuité des
  racines d'un polynôme à racines simples).

  Soient $r,r'$ deux réflexions de $\NCP_W$ conjuguées. Par surjectivité de
  $\fact$, il existe $y,y' \in Y -\CK$, tels que
  $\fact(y)=(r,r_2,\dots,r_n)$ et $\fact(y')=(r',r_2',\dots,r_n')$. On peut
  supposer que $\LL(y)=\LL(y')$ ; soit $x$ leur élément minimal pour $\lex$
  (correspondant à $r$ et $r'$). Soit $\Lambda$ la strate de $\Lb_1$
  correspondant à la classe de conjugaison de $r$ et $r'$. D'après le lemme
  \ref{lemtype}, $(y,x)$ et $(y',x)$ sont dans $\Lambda ^0$. Plus
  précisément, $(y,x)$ et $(y',x)$ sont dans $\Lambda \cap \CH'$, que l'on
  va noter $\Lambda'$.

  Notons que $\Lambda'$ est un ouvert de Zariski de $\Lambda$, donc est
  connexe par arcs, par le même argument que pour la proposition
  \ref{propcpa}. Par conséquent, on peut relier $(y,x)$ et $(y',x)$ par un
  chemin $\gamma$ dans $\CH '$. Celui-ci se projette en un chemin dans $Y-
  \CK$, et détermine via $\LL$ un lacet dans $\Enreg$. Ce lacet
  représente une tresse $\beta$ qui, par construction, stabilise le premier
  brin $(x)$. Ainsi $y'=y \cdot \beta$ par l'action de monodromie, et
  $\fact(y')=\fact(y)\cdot \beta$ par l'action d'Hurwitz. Comme $\beta$
  stabilise le brin $(x)$, on en déduit que $r$ et $r'$ sont fortement
  conjugués, en vertu de la remarque \ref{rqrefl}.
\end{proof}


\begin{thebibliography}{99}

\bibitem{Armstrong} D. Armstrong, \emph{Generalized noncrossing partitions
    and combinatorics of Coxeter groups}, Thesis (2006),
  Mem. Amer. Math. Soc. \textbf{202} (2009), no. 949.

\bibitem{arnold} V. I. Arnold, \emph{Critical points of functions and
    the classification of caustics}, Russian Math. Surveys \textbf{29}
  (1974), 243--244.

\bibitem{athareiner} Ch. Athanasiadis, V. Reiner,
\emph{Noncrossing partitions for the group $D_n$},
SIAM J. Discrete Math. \textbf{18} (2004), 397--417.

\bibitem{zariski} D. Bessis, \emph{Zariski theorems and diagrams for
    braid groups}, Invent. Math. \textbf{145} (2001), 487--507.

\bibitem{Bessisdual} \bysame, \emph{The dual braid monoid}, Ann.
  Sci. Éc. Norm. Supér. (4) \textbf{36} (2003), 647--683.

\bibitem{BessisKPi1} \bysame, \emph{Finite complex reflection
    arrangements are $K(\pi,1)$}, arXiv preprint
  \href{http://arxiv.org/abs/math/0610777}{\texttt{math.GT/0610777}}.


\bibitem{BessisCorran} D. Bessis, R. Corran, \emph{Non-crossing
    partitions of type $(e,e,r)$}, Adv. Math. \textbf{202} (2006), 1--49.

\bibitem{BDM} D. Bessis, F. Digne, J. Michel, \emph{Springer theory in
    braid groups and the Birman-Ko-Lee monoid}, Pacific J. Math. \textbf{205}
  (2002), 287--310.

\bibitem{BKL} J. Birman, K. H. Ko, S. J. Lee, \emph{A new approach to
    the word and conjugacy problem in the braid groups}, Adv. Math. \textbf{
    139} (1998), no. 2, 322--353.

\bibitem{BWordre} T. Brady, C. Watt, \emph{A partial order
    on the orthogonal group}, Comm. Algebra \textbf{30} (2002), no.
  8, 3749--3754.

\bibitem{BWordre2} \bysame, \emph{$K(\pi,1)$'s for
    Artin groups of finite type}, Proceedings of the
  Conference on Geometric and Combinatorial Group Theory,
  Part I (Haifa, 2000), Geom. Dedicata \textbf{94} (2002), 225--250.

\bibitem{BWtreillis} \bysame, \emph{Non-crossing partition
    lattices in finite real reflection groups},
  Trans. Amer. Math. Soc. \textbf{360} (2008), no. 4, 1983--2005. 

\bibitem{Brieskorn} E. Brieskorn, \emph{Automorphic Sets and Braids and
    Singularities}, Contemp. Math. \textbf{78}, 1988 (Braids), 45--115.


\bibitem{broumarou1} M. Broué, G. Malle, R. Rouquier, \emph{On
    complex reflection groups and their associated braid groups},
  Representations of groups (Banff, AB, 1994), CMS Conf. Proc., vol. 16,
  Amer. Math. Soc., Providence (1995), 1--13.

\bibitem{broumarou} \bysame, \emph{Complex
    reflection groups, braid groups, Hecke algebras}, J. Reine
    Angew. Math. \textbf{500} (1998), 127--190.

\bibitem{carter} R. W. Carter, \emph{Conjugacy classes in the Weyl
  groups}, Compos. Math. \textbf{25} (1972), 1--52.

\bibitem{chapoton} F. Chapoton, \emph{Enumerative properties of
    generalized associahedra}, Sém. Lothar. Combin.  \textbf{51} (2004),
  Art. B51b (electronic).

\bibitem{LD} P. Dehornoy, \emph{Braids and self-distributivity},
  Progress in Math. \textbf{192}, Birkhäuser, 2000.

\bibitem{Deh2} \bysame, \emph{Groupes de Garside}, Ann. Sci. Éc.
  Norm. Supér. (4) \textbf{35} (2002), 267--306.

\bibitem{deligneletter} P.~Deligne, letter to E.~Looijenga, 9/3/1974.

\bibitem{IngThom} C. Ingalls and H. Thomas, \emph{Noncrossing partitions
    and representations of quivers}, arXiv preprint
  \href{http://arxiv.org/abs/math.RT/0612219}{\texttt{math.RT/0612219}}, to
  appear in Compos. Math. .

\bibitem{LZgraphs} S. K. Lando, A. K. Zvonkin, \emph{Graphs on surfaces
    and their applications}, Encyclopaedia Math. Sci. \textbf{141},
  Springer, 2004.

\bibitem{landozvonkine} S. K. Lando, D. Zvonkine, \emph{On
    multiplicities of the Lyashko-Looijenga mapping on the discriminant
    strata}, Funct. Anal. Appl. \textbf{33} (1999), no.
  3, 21--34 and 96 (in Russian), 178--188 (in English).

\bibitem{looijenga} E. Looijenga, \emph{The complement of the bifurcation
    variety of a simple singularity}, Invent. Math. \textbf{23} (1974),
  105--116.

\bibitem{reading} N. Reading, \emph{Chains in the noncrossing partition
    lattice}, SIAM J. Discrete Math. \textbf{22} (2008), no. 3,
  875--886. 

\bibitem{reiner} V. Reiner, \emph{Non-crossing partitions for classical
    reflection groups}, Discrete Math. \textbf{177} (1997), 195--222.

\bibitem{Springer} T. A. Springer, \emph{Regular elements of finite
    reflection groups}, Invent. Math. \textbf{25} (1974), 159--198.

\bibitem{stanley} R. P. Stanley, \emph{Enumerative combinatorics}, vol.
  1, Cambridge Stud. Adv. Math. \textbf{29}, Cambridge Univ. Press, 1997.

\bibitem{steinberg} R. Steinberg, \emph{Differential equations
    invariant under finite reflection groups}, Trans. Amer. Math.
  Soc. \textbf{112} (1964), 392-–400.

\end{thebibliography}
\end{document}